\documentclass[letterpaper, 10 pt, conference]{ieeeconf}
\IEEEoverridecommandlockouts

\usepackage{amsmath}
\usepackage{graphicx}
\usepackage{cite}
\usepackage{authblk}
\usepackage{amssymb}
\usepackage{subfigure}
\usepackage{multirow}
\usepackage{color}	
\usepackage{subfigure}
\usepackage{algorithmic,algorithm}
\usepackage{gensymb}
\usepackage{tikz} \usetikzlibrary{arrows,calc,decorations.pathreplacing}

\newtheorem{theorem}{Theorem}

\newtheorem{corollary}{Corollary}

\title{The Complete Differential Game of Active Target Defense}
\author{Eloy Garcia,  David W. Casbeer, and Meir Pachter
\thanks{E. Garcia and D. Casbeer are with the Control Science Center of Excellence, Air Force Research Laboratory, Wright-Patterson AFB, OH 45433. Corresponding author \ttfamily{eloy.garcia.2@us.af.mil}}
\thanks{M. Pachter is with the Department of Electrical Engineering, Air Force Institute of Technology, Wright-Patterson AFB, OH 45433.}
}

\begin{document}
\maketitle 
\begin{abstract}
In the Target-Attacker-Defender (TAD) differential game, an Attacker missile strives to capture a Target aircraft. The Target tries to escape the Attacker and is aided by a Defender missile which aims at intercepting the Attacker before the latter manages to close in on the Target. The conflict between these intelligent adversaries has been suitably modeled as a zero-sum differential game. Optimal strategies have been synthesized covering the region of the state space where the Target/Defender team is able to win the game. However, the Game of Degree in the Attacker's region of win has not been fully addressed. 
Preliminary attempts at designing the players' strategies have not been proven to be optimal in the differential game sense. 
The main results of the paper present the optimal strategies of the Game of Degree in the Attacker's winning region of the state space. It is proven that the obtained strategies provide the saddle-point solution of the game; the Value function is obtained and it is shown to be continuous and continuously differentiable. It is also demonstrated that it is the solution of the Hamilton-Jacobi-Isaacs (HJI) equation. 
Finally, the obtained strategies are compared to recent results addressing the TAD differential game in  \cite{Liang19}. It is shown by counterexample that the strategies proposed in \cite{Liang19} are not optimal. The unique regular solution of this differential game that actually provides a semipermeable Barrier surface is synthesized and verified in this paper. 
\end{abstract}

\section{Introduction} \label{sec:intro}
In conflict and combat scenarios involving autonomous agents, the synthesis of intelligent actions must consider the potential strategies by the adversary. Uncertainty on the adversary's behavior, actions, and/or decisions represents one of the main challenges when analyzing possible outcomes of conflict situations.
Differential game theory \cite{Isaacs65,Basar99} provides the right framework to analyze dynamic conflict and combat scenarios and to design optimal strategies for each one of the players.

Differential game theory has been applied to the study of pursuit-evasion problems  \cite{Huang11}, \cite{Oyler16}. 
Pursuit-evasion problems also arise in Unmanned Aerial Vehicles (UAVs) operations such as in \cite{Sprinkle04} where a receding-horizon approach provides evasive maneuvers for a UAV assuming a known pursuer's input. In \cite{EarlDandrea07}, a multi-agent scenario is considered where a number of pursuers are tasked to intercept a group of evaders and where the goals of the evaders are assumed to be known. 

The interesting problem in \cite{breakwell1979point} considered the dynamic game of a fast pursuer trying to capture in minimal time two slower evaders in succession. The evaders, on their part, cooperate and try to maximize the capture time. This work was extended in \cite{liu2013evasion} where the fast pursuer attempts to sequentially capture several evaders.
Similarly, the slow evaders act as a team and cooperate in order to maximize the total time from the beginning of the game until the last evader is captured. The numerical solution provided in that reference shows that the optimal strategies of every player consist of constant headings (the pursuer's heading is piecewise constant as one could expect and it changes at time instants when an evader is captured). 

The problem of active defense of a target aircraft is studied in this paper by means of differential game theory. The Target ($T$), the Attacker ($A$), and the Defender ($D$) are the participants. In this scenario, $A$ pursues and tries to capture $T$. $D$ aims at intercepting $A$ before the latter can reach $T$. Hence, agents $D$ and $T$ cooperate and form a team to achieve interception of $A$ by $D$ and successful escape of $T$. The Target-Attacker-Defender (TAD) differential game is a two-termination set differential game where two distinct outcomes can be obtained. The first outcome is given by interception of $A$ by $D$ and $T$ escapes. In this case the $T/D$ team wins the game. The second outcome is given by $A$ capturing $T$ before being intercepted by $D$. In such a case the Attacker wins the game. These two outcomes give rise to the Game of Kind and its solution partitions the state space into two regions of win: $\mathcal{R}_e$, the escape region, and $\mathcal{R}_c$, the capture region. This paper focuses on synthesizing the optimal strategies of the players in the region $\mathcal{R}_c$, where the Attacker wins the game under optimal play.

Several references have considered the active defense of aircraft. For instance, Li and Cruz \cite{LiCruz11} considered the game of defending an asset from an attacking intruder using an interceptor.  
The differential game of active target defense in the presence of obstacles was analyzed in \cite{fisac2015pursuit} and the case where the defender is endowed with a positive capture radius was addressed in \cite{Garcia18TAES}. 
Active defense when the Target aircraft is non-maneuverable was investigated in \cite{harini2015new} and \cite{Weintraub18}. 
A differential game with multiple attackers, multiple defenders, and a stationary target in a bounded domain was analyzed in \cite{chen2016multiplayer}. Due to numerical intractability the authors of this reference use the solution of the single attacker-single defender case in order to determine pairwise outcomes favorable to the defender team.
The work in \cite{zhou2016cooperative} considers a group of cooperative pursuers that try to capture a single evader within a bounded domain. 
The papers \cite{Lorenzetti18,zhou2012,margellos2011hamilton,zhou2018efficient} provide different approaches to solve reach and avoid games.

The Active Target Defense Differential Game (ATDDG) was introduced in \cite{Pachter14} and further analyzed in \cite{Garcia15JGCD}; it solves the TAD differential game in $\mathcal{R}_e$, the $T/D$ team's winning region. On the other hand the TAD differential game in the region $\mathcal{R}_c$, the Attacker's winning region, was addressed in the recent Reference \cite{Liang19}. However, that reference proposed an incorrect solution. The authors of \cite{Liang19} failed to provide the important property of a semipermeable surface. In other words, the strategies designed in \cite{Liang19} do not guarantee that the Attacker wins when the state of the system is in $\mathcal{R}_c$, where the Attacker is supposed to win. 

We show by counterexample that when the state of the system is in the Attacker's winning region $\mathcal{R}_c$ and the Attacker implements the strategy proposed in \cite{Liang19}, then, there exists a non-anticipative  strategy for the $T/D$ team to actually defeat the Attacker and win the game.
The $T/D$ team implements the optimal strategy derived in this paper and the Attacker loses the game since it \textit{does not} capture the Target. However, if the Attacker employs the strategy obtained in this paper then it can capture the Target regardless of the strategy implemented  by the $T/D$ team; the best the Target and Defender can do is to implement the solution in this paper which is the saddle-point solution of the TAD differential game.

The papers \cite{Garcia18ACC,Garcia18pursuit} proposed a set of optimal strategies for the TAD differential game in the Attacker's winning region employing a reduced state space and using the geometric properties of the problem. However, proof of optimality was absent in these papers; this issue is addressed in the current paper by finding the Value function, the solution of the  Hamilton-Jacobi-Isaacs (HJI) Partial Differential Equation (PDE).  
The possible non-differentiability of the solution of the HJI PDE is a concern in differential games and a generalized solution concept is provided by the viscosity solution \cite{Crandall84}, \cite{Lions85}. 
However, it is important to note that non-differentiability is not an issue in this paper.  By using Isaacs' method \cite{Isaacs65} the players' optimal strategies are synthesized. Absence of singular surfaces is demonstrated and the Value function is $C^1$. As it has been stated in \cite{chen2016multiplayer}, Isaacs' method is the ideal situation in differential games, if it is attainable. Unfortunately, Isaacs' method does not scale well as the dimension of the state increases. Many games may still have a classical solution in closed-form but it is very difficult to obtain it. Hence, when the state space of the differential game is of higher dimension, obtaining a closed-form $C^1$  Value function which is the classical solution of the HJI PDE represents a valuable contribution.

The main contributions of this paper are as follows. A complete treatment of the TAD differential game in the Attacker's winning region is presented. Here, we provide a rigorous synthesis of the players' optimal strategies by means of Pontryagin's Maximum Principle as opposed to the geometric approach utilized in \cite{Garcia18ACC,Garcia18pursuit}. Additionally, we prove that the obtained solution is indeed the solution of the differential game under consideration by obtaining the Value function and showing that the Value function is continuous, continuously differentiable, and it is the solution of the HJI PDE. Furthermore, we show that the obtained solution is equivalent to the ATDDG optimal solution when the state of the system is located on the Barrier surface that separates the two winning regions. Finally, we compare the obtained strategies in this paper to recent results obtained in \cite{Liang19}. Both, in that reference and in this paper, two outcomes of the game are defined: 1) $T$ is captured and $A$ wins; 2) $A$ is intercepted by $D$ before reaching $T$ and the $T/D$ team wins. The Game of Kind is similarly solved in terms of terminal distance between the players. However, capture time is used as the objective function in \cite{Liang19} which results in Pure Pursuit strategy for the Attacker to implement in its winning region. Pure Pursuit  \textit{is not} the optimal strategy when the Attacker pursues the Target in the presence of the Defender since we are able to show that when the state of the game is initially in the Attacker's winning region, there exist an strategy (that we also derived in this paper) such that the Defender intercepts the Attacker and the Target escapes; the $T/D$ team wins the game.
Also important is the fact that the strategies in \cite{Liang19} do not result in a semipermeable barrier surface as it was assumed in that reference. This is of great importance since one of the teams is capable of switching the regions (make the state to cross the barrier surface and change the outcome of the game) by applying the optimal strategy in this paper while the opponent applies the strategy in \cite{Liang19}.

The paper is organized as follows. Section \ref{sec:Problem} states the Target-Attacker-Defender (TAD) Differential Game and it provides the solution to the Game of Kind. The Game of Degree in the Attacker's winning region is solved in Section \ref{sec:num}. The properties of the game on the Barrier surface are analyzed in Section \ref{sec:Barrier}. Illustrative examples are shown in Section \ref{sec:examples}. Concluding remarks are included in Section \ref{sec:concl}.

\section{The TAD Differential Game} \label{sec:Problem}

\subsection{Two-termination set differential game}
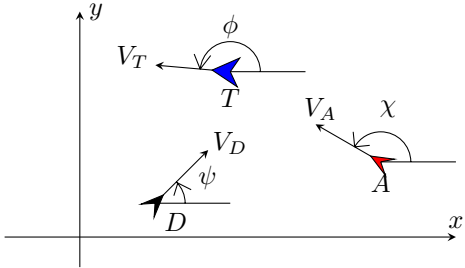
\begin{figure}[htb]\centering
\begin{tikzpicture}[>=stealth]
	\tikzset{
	   uav/.pic = {  
	      \draw [pic actions] (0,0)
	         -- (-1mm,1.75mm)
	         -- ++(3.2mm,-1.75mm)  
	         -- ++(-3.2mm,-1.75mm)
	         -- (0,0);
	   }
         }
	\tikzset{
	   missile/.pic = {  
	      \draw [pic actions] (0,0)
	         -- (-2mm,1.75mm)
	         -- ++(4.2mm,-1.75mm)  
	         -- ++(-4.2mm,-1.75mm)
	         -- (0,0);
	   }
	}
	\coordinate (D) at (1,.45);  
	\coordinate (T) at (2,2.2);
	\coordinate (A) at (4,1);
	\draw [->] (D) -- ++ (0:1)
			(D) ++(0:.4) arc [start angle=0, end angle=45, radius=.4cm]
	    (1.28,.73) -- ++(.2,-.1)
  		(1.3,.71) -- ++(-.05,-.16)
		        (1.7,.8) node{$\psi$}
						(2,1.25) node{$V_D$}
		        (D) -- ++(45:1);
	\draw (D) pic [rotate=45,fill=black,scale=.7] {missile};  
	\draw (D) node [below right] {$D$};
	\draw [->] (A) -- ++ (0:1)
			(A) ++(0:.4) arc [start angle=0, end angle=-30+180, radius=.4cm]
		  (3.654,1.2) -- ++(-.02,.2)
			(3.654,1.2) -- ++(.2,.02)
						(4.1,1.7) node{$\chi$}
						(3.2,1.7) node{$V_A$}
		        (A) -- ++(-30+180:1);
	\draw (A) pic [rotate=-30+180,fill=red,scale=.7] {missile};
	\draw (A) node [below=1] {$A$};
	\draw [->] (T) -- ++ (0:1)
			(T) ++(0:.4) arc [start angle=0, end angle=175, radius=.4cm]
		  (1.60,2.235) -- ++(-.08,.2)
			(1.60,2.235) -- ++(.14,.18)
			      (2,2.8) node{$\phi$}
						(.7,2.4) node{$V_T$}
		        (T) -- ++(175:1);
	\draw (T) pic [rotate=175,fill=blue,scale=1.1] {uav};
	\draw (T) node [below=3] {$T$};
	\node [right] (Y) at (0,3) {$y$};
	\draw [->](0,-.4) -- (0,3);
	\node [above] (X) at (5,0) {$x$};
	\draw [->](-1,0) -- (5,0);
 \end{tikzpicture}
\caption{Target, Attacker, and Defender scenario}
\label{fig:problem description}
\end{figure}

The scenario of active target defense considers three players: a Target ($T$), an Attacker ($A$), and a Defender ($D$) which have ``simple motion" $\grave{\text{a}}$ la Isaacs, they are holonomic.  The game is played in the Euclidean plane where the controls of $T$, $A$, and $D$ are their respective instantaneous headings $\phi$, $\chi$, and $\psi$ and their states are specified by their Cartesian coordinates $\textbf{x}_T=(x_T,y_T)$, $\textbf{x}_A=(x_A,y_A)$, and $\textbf{x}_D=(x_D,y_D)$, as it is shown in Fig. \ref{fig:problem description}. The players $T$, $A$ and $D$ have constant speeds denoted by $V_T$, $V_A$, and $V_D$, respectively.  The complete state of the TAD differential game is specified by $\textbf{x}:=( x_T, y_T, x_A, y_A, x_D, y_D)\in \mathbb{R}^6$. The game set is the entire space $\mathbb{R}^6$. The initial time is denoted by $t_0$ and the corresponding initial state of the system is 
\begin{align}
	\textbf{x}_0 := (x_{T_0}, y_{T_0}, x_{A_0}, y_{A_0}, x_{D_0}, y_{D_0}) = \textbf{x}(t_0).  \nonumber
\end{align}
The Target aircraft is slower than the Attacker missile, then the speed ratio $\alpha=V_T/V_A<1$. When $V_A<V_T$ there is no need for a Defender, the Target always outruns the Attacker. 
We assume that the Attacker and Defender have similar capabilities, so $V_A=V_D$. Without loss of generality, the players' speeds are normalized so that $V_A=V_D=1$ and $V_T=\alpha$.

\textit{Remark}. When $V_A>V_D$ (slower Defender) and point capture is required then $A$ always captures $T$ irrespective of the initial conditions. The slower agent, $D$ in this case, is incapable of achieving point interception of the faster agent, $A$. Player $A$ can always exploit its speed advantage to circumvent a slowly moving point and capture $T$. Additionally, since $\alpha<1$ ($T$ is slower than $A$) the $T/D$ team is not able to indefinitely keep $A$ away from $T$ or, in terms of \cite{Liang19}, keep the state such that $R>0$ and $r>0$ where $R$ is the $A-T$ separation and $r$ is the $A-D$ separation. Hence, a rendezvous strategy between $T$ and $D$ always results in $T$ being captured by $A$. Therefore, in this paper we focus on the case $V_A=V_D$. The results can also be extended to the case $V_D>V_A$.

The control input of the $T/D$ team is the pair of instantaneous headings $\textbf{u}_{T,D}=\left\{\phi,\psi\right\}$.
The Attacker's control is his instantaneous heading angle, $\textbf{u}_A=\left\{\chi\right\}$.  The dynamics $\dot{\textbf{x}}=\textbf{f}(\textbf{x},\textbf{u}_A,\textbf{u}_{T,D})$ are specified by the system of ordinary differential equations
\begin{align}
 \left.
	\begin{array}{l l}
        \dot{x}_A&=\cos\chi,    \ \ \ \ \ x_A(0)=x_{A_0}   \\
	\dot{y}_A&=\sin\chi,  \ \ \ \ \ y_A(0)=y_{A_0}  \\
	\dot{x}_D&=\cos\psi,   \ \ \ \ \ x_D(0)=x_{D_0} \\
	\dot{y}_D&=\sin\psi,   \ \ \ \ \ y_D(0)=y_{D_0}  \\
	\dot{x}_T&=\alpha\cos\phi,  \ \ \  x_T(0)=x_{T_0}      \\
	\dot{y}_T&=\alpha\sin\phi,    \ \ \  y_T(0)=y_{T_0}    
	\end{array}  \right.   \label{eq:xT}
\end{align}
where $\alpha<1$ is the problem parameter and the admissible controls are given by $\chi,\phi,\psi \in [-\pi,\pi]$. Both, the state and the controls, are unconstrained.

The terminal condition is point capture, that is, the separation between Target and Attacker becomes zero allowing the Attacker to capture the Target and win the game. An alternative termination condition is  when the separation between Attacker and Defender is equal to zero; this case represents interception of the Attacker by the Defender  and the $T/D$ team wins. Hence, the termination set for the complete TAD differential game is
\begin{align}
   \mathcal{S} :=  \mathcal{S}_e    \   \bigcup \   \mathcal{S}_c   \label{eq:TwoSets}
\end{align}
where 
\begin{align}
   \mathcal{S}_e:=  \big\{ \ \textbf{x} \ | \sqrt{ (x_A-x_D)^2 + (y_A-y_D)^2} =0  \big\}   \nonumber
\end{align}
represents interception of the Attacker by the Defender (and the Target escapes) and 
\begin{align}
   \mathcal{S}_c:=   \big\{ \ \textbf{x} \ |  \sqrt{ (x_A-x_T)^2 + (y_A-y_T)^2}=0 \big\}    \nonumber
\end{align}
represents the opposite outcome where the Attacker wins by capturing the Target. Conceptually, the TAD differential game with termination set given in \eqref{eq:TwoSets} belongs to the class of two-target or two-termination set differential games which have been considered early by different authors \cite{getz1979qualitative,Getz1981capturability,ardema1985combat}. This two-termination set differential game concept was introduced in order to extend classical pursuit-evasion games where only one target, or termination, set exists. Naturally, the two-target differential game is useful in the analysis of combat games such as \cite{ardema1985combat} and also the differential game under consideration where the roles of pursuer and evader are not designated ahead of time; instead each player wants to defeat his opponent by terminating the game in its own terminal set. In the TAD differential game the Attacker pursues the Target but it is also being pursued by the Defender. These roles of $A$ and the associated termination sets \eqref{eq:TwoSets} suitably describe the possible outcomes of the TAD scenario.

\subsection{Game of Kind}  \label{subsec:gok}
Due to the two different outcomes of the TAD differential game specified in  \eqref{eq:TwoSets}, the Game of Kind needs to be addressed and solved in order to partition the state space into two winning regions, one for each team. Since different Games of Degree are played in each region, it is essential for each player to determine which region the current state of the system lies in, then the corresponding Game of Degree is solved and the appropriate optimal strategies are implemented.
 The notions of Game of Degree and Game of Kind are central in pursuit and evasion games \cite{Isaacs65}. The state space $\mathbb{R}^6$ is partitioned into two sets: $\mathcal{R}_e$ and $\mathcal{R}_{c}$. Let us note first that when the Target is closer to the Defender than to the Attacker, the Target can always escape; even in the case when $\alpha \rightarrow 0$ since the Target is already located in the dominance region of the Defender. Thus, the closed half-space $\mathcal{R}_{ed}$ belongs to the escape set $\mathcal{R}_e$, that is, $\mathcal{R}_{ed} \subset  \mathcal{R}_e$, where 
 \begin{align}
\left.
	 \begin{array}{l l}
\mathcal{R}_{ed}:= & \big\{ \ \textbf{x} \ | \sqrt{ (x_A-x_T)^2 + (y_A-y_T)^2} \\
&\qquad -\sqrt{ (x_D-x_T)^2 + (y_D-y_T)^2} \geq 0  \big\}. 
\end{array}   \right.  \label{eq:Red}
\end{align}

 
  We define the sets $\mathcal{R}_{ea}$ and $\mathcal{R}_c$ as follows 
\begin{align}
\left.
	 \begin{array}{l l}
\mathcal{R}_{ea}:=  \big\{ \ \textbf{x} \ | \textbf{x}\notin \mathcal{R}_{ed} , \ B( \textbf{x};\alpha)<0  \big\} \\
\mathcal{R}_c:=  \big\{ \ \textbf{x} \ | \textbf{x}\notin \mathcal{R}_{ed} , \  B( \textbf{x};\alpha)>0  \big\}  
\end{array}   \right.  \label {eq:ReDef}
\end{align}
and the escape set is then 
\begin{align}
\left.
	 \begin{array}{l l}
\mathcal{R}_e:=  \mathcal{R}_{ed} \cup \mathcal{R}_{ea}.
\end{array}   \right.  \label {eq:EscapeSet}
\end{align}
The Barrier surface is defined as
\begin{align}
\left.
	 \begin{array}{l l}
\mathcal{B}:=  \big\{ \ \textbf{x} \ | \textbf{x}\notin \mathcal{R}_{ed} , \ B( \textbf{x};\alpha)=0  \big\} \\
\end{array}   \right.  \label {eq:BarrSurface}
\end{align}
where the Barrier function $B( \textbf{x};\alpha)$ is explicitly characterized as follows.
%

\begin{theorem}
For a given speed ratio parameter $0<\alpha<1$, the Barrier surface that separates the state space $\mathbb{R}^6$ into the two regions $\mathcal{R}_e$ and $\mathcal{R}_{c}$ is given by $B( \textbf{x};\alpha)=0$, where
\begin{align}
\left.
	 \begin{array}{l l} 
	B( \textbf{x};\alpha)&= b_{xx} x_T^2 + b_{yy}y_T^2+2b_{xy}x_Ty_T \\
	&~~  +2b_x x_T+2b_y y_T+ b.  
	\end{array}   \right.   \label{eq:hb}
\end{align}
The coefficients of \eqref{eq:hb} are given by
\begin{align}
\left.
	 \begin{array}{l l} 
   b_{xx}= \cos^2\sigma-\alpha^2 \\
	 b_{yy}= \sin^2\sigma-\alpha^2  \\
	 b_{xy}= \sin\sigma\cos\sigma \\
	 b_x= \alpha^2x_A\sin^2\sigma - (1-\alpha^2)x_0\cos^2\sigma \\
	 \quad \ \ \  -[\alpha^2y_A +(1-\alpha^2)y_0]\sin\sigma\cos\sigma \\
	 b_y= \alpha^2y_A\cos^2\sigma - (1-\alpha^2)y_0\sin^2\sigma \\
	  \quad \ \ \  -[\alpha^2x_A +(1-\alpha^2)x_0]\sin\sigma\cos\sigma \\
	 b= \big( [\alpha^2y_A +(1-\alpha^2)y_0]\sin\sigma \\
	 \quad \ \ + [\alpha^2x_A +(1-\alpha^2)x_0]\cos\sigma \big)^2 -\alpha^2(x_A^2\!+\!y_A^2) 
\end{array}   \right.  \label{eq:hbcoef}
\end{align}
where 
\begin{align}
\left.
	 \begin{array}{l l} 
 \cos\sigma=\frac{x_A-x_D}{\sqrt{(x_A-x_D)^2+(y_A-y_D)^2}} \\
 \sin\sigma=\frac{y_A-y_D}{\sqrt{(x_A-x_D)^2+(y_A-y_D)^2}} \\
 x_0(\textbf{x}) = \frac{1}{2}(x_A+x_D) \\
 y_0(\textbf{x}) = \frac{1}{2}(y_A+y_D).
 \end{array}   \right.  \nonumber
\end{align}
\end{theorem}
For proof and additional details see \cite{Garcia19TAC}.

\begin{figure}
	\begin{center}
		\includegraphics[width=8.0cm,trim=1.0cm .0cm 1.0cm .0cm]{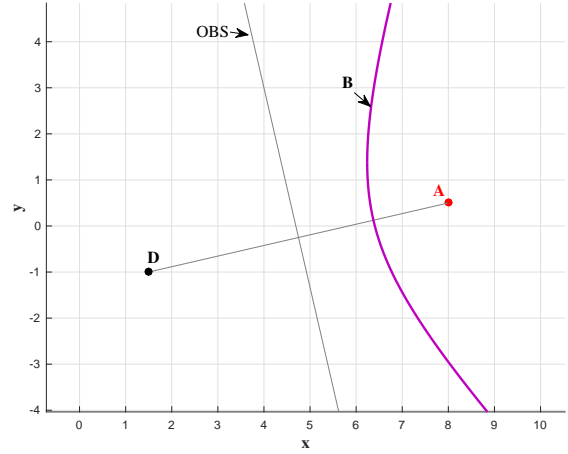}
	\caption{Cross section of the Barrier surface }
	\label{fig:BEX}
	\end{center}
\end{figure}

\textit{Remark}. The choice of notation in \eqref{eq:hb} is due to the fact that the cross section of the Barrier surface $\mathcal{B}$, for fixed $A$ and $D$ positions, represents an hyperbola. This provides a clear illustration of the points in the Cartesian plane for which the Target is guaranteed to escape under optimal play -- see Fig. \ref{fig:BEX}. Note that only the $A$ branch of the hyperbola is relevant in the solution to the Game of Kind. The other branch is irrelevant since it is located in $D$'s dominance region where $T$ can always escape since $\textbf{x} \in \mathcal{R}_{ed}$ holds. $D$'s dominance region in Fig. \ref{fig:BEX} is delineated by the orthogonal bisector (OBS) of the segment $A-D$.

If the state is such that $\textbf{x} \in \mathcal{R}_e$ then the ATDDG Game of Degree is played where the opposing teams try to min/max the terminal $A-T$ separation; this Game of Degree was formulated and solved in \cite{Garcia19TAC}. 
By contrast, in this paper we focus on the Game of Degree in the Attacker's winning region where the state is such that  $\textbf{x} \in \mathcal{R}_c$. We refer to the Game of Degree in the  Attacker's winning region as the Capture Differential Game (CDG). In the CDG, the opposing teams try to min/max the terminal $D-T$ separation.

\subsection{The Game of Degree in $\mathcal{R}_c$}

The CDG is a perfect information game where every player knows the dynamics \eqref{eq:xT} and the speed ratio parameter $\alpha$. The players have access to the state $\textbf{x}$ at the current time $t$. The optimal strategies will be state feedback strategies. Finally, and most importantly, it is assumed that the agents do not know the opponent's current decision: no discriminatory/stroboscopic strategies are needed in the CDG.

The termination condition in the CDG is
\begin{align}
 \left.
	\begin{array}{l l}
   x_A=x_T,  \ \ \ \ \     y_A=y_T.	
   \end{array}  \right. 
\label{eq:PCx} 
\end{align}
The terminal manifold \eqref{eq:PCx} is a four dimensional hyperplane in $\mathbb{R}^6$.
The terminal time $t_f$ is \textit{not} fixed and is instead defined as the time instant when the state of the system satisfies \eqref{eq:PCx}. At time instant $t_f$ the terminal state is $\textbf{x}_f: = (x_{T_f}, y_{T_f}, x_{A_f}, y_{A_f}, x_{D_f}, y_{D_f}) = \textbf{x}(t_f)$.
The terminal cost/payoff is
\begin{align}
  J(\textbf{u}_A(t),\textbf{u}_{T,D}(t);\textbf{x}_0)=\Phi(\textbf{x}_f)	\label{eq:costDG2}
\end{align}
where
\begin{align}
  \Phi(\textbf{x}_f):=\sqrt{(x_{D_f}-x_{T_f})^2+(y_{D_f}-y_{T_f})^2}.	\label{eq:costDG}
\end{align}
The cost/payoff functional depends only on the terminal state; its Value is given by
\begin{align}
  V(\textbf{x}_0):= \max_{\textbf{u}_A(\cdot)} \ \min_{\textbf{u}_{T,D}(\cdot)} J(\textbf{u}_A(\cdot),\textbf{u}_{T,D}(\cdot);\textbf{x}_0)	\label{eq:costDG3}
\end{align}
subject to \eqref{eq:xT} and \eqref{eq:PCx}, where $\textbf{u}_A(\cdot)$ and $\textbf{u}_{T,D}(\cdot)$ are the teams' state feedback strategies.
In detail, the Attacker aims at capturing the Target and maximizing the terminal $T-D$ separation. The Target and Defender, knowing that capture of the Target is imminent under optimal play, cooperate in order to minimize the terminal $T-D$ separation. 

The synthesis of optimal strategies of the CDG, the Game of Degree when $\mathbf{x} \in \mathcal{R}_c$, is important because, together with the solution of the ATDDG in \cite{Garcia19TAC,Pachter16}, the Game of Degree when  $\mathbf{x} \in \mathcal{R}_e$, they provide the three players' optimal strategies everywhere in the state space $\mathbb{R}^6=\mathcal{R}_c \ \cup \ \mathcal{R}_e$.  Therefore, the complete solution of the three-agent pursuit-evasion TAD differential game which comprises a Target, an Attacker, and a Defender is obtained. 

The cost/payoff functional  \eqref{eq:costDG2} is designed to mesh with the cost/payoff functional defined in \cite{Garcia19TAC}.  The ATDDG is played in that case since $\textbf{x} \in \mathcal{R}_e$. The solutions of both Games of Degree, the ATDDG in \cite{Garcia19TAC} and the CDG in this paper, possess important properties related to two-termination set differential games. The complete optimal strategies comprising both Games of Degree result in a semipermeable Barrier surface $\mathcal{B}$ that separates the two winning regions. This semipermeable surface property is not obtained in  \cite{Liang19} where the same problem is considered. This will be discussed in Section \ref{sec:examples}.

The complete solution of the TAD differential game (consisting of the solutions of the ATDDG and the CDG) yields state feedback strategies for the players which possess the continuity property expected in two termination set differential games. In more detail,
when the state of the system $\textbf{x}$ is located on the Barrier surface, both, the CDG and the ATDDG, have feasible solutions; furthermore, simultaneous capture of $T$ by $A$ and interception of $A$ by $D$ is attained at termination. Given the nature of the terminal cost functionals in the CDG and in the ATDDG, the solution of the Game of Degree in $\mathcal{R}_c$ needs to provide the same strategies and the same Value of the game as the solution of the Game of Degree in $\mathcal{R}_e$ when simultaneous capture is realized, \textit{i.e.} $A$ captures $T$ at the same time instant when $D$ intercepts $A$. Thus, the game terminates on the subspace $ \mathcal{S}_e   \bigcap    \mathcal{S}_c$. This  will be corroborated in Section \ref{sec:Barrier}.

\textit{Remark}. The choice of cost/payoff in the CDG has also a practical interpretation. Since it is specified by the solution of the Game of Kind  that $A$ will capture $T$ despite the best efforts by the $T/D$ team, then the logical strategy for $D$ is to try to reach a point at time $t_f$ which is as close as possible to $T$. In the case $A$ makes a mistake, the strategy of $T\ \& \ D$ cooperating to minimize their terminal separation provides a reasonable strategy and immediately takes advantage of non-optimal behaviors of $A$. If player $A$ pursues $T$ using a different guidance law other than the optimal strategy derived in this paper, such as the Pure Pursuit guidance proposed in \cite{Liang19}, then, the terminal $T-D$ distance will decrease with respect to the guaranteed Value of the game, which is attained when all players act optimally. 
And under further non-optimal play by the Attacker, the state can cross the Barrier surface and it will then hold that $\textbf{x}\in \mathcal{R}_e$ where the strategies of the ATDDG, the Game of Degree in $\mathcal{R}_e$, will be employed by $T$ and $D$, and now the Target can escape.

\section{Optimal Strategies in $\mathcal{R}_c$}  \label{sec:num}
In this section we provide the players' optimal strategies in the Attacker's winning region. The solution is synthesized by employing the  ``Two-sided'' Pontryagin Maximum Principle (PMP) which is also known as the ``Two-person'' extension of the PMP \cite{Basar99}.  This is justified due to the absence of singular surfaces in the CDG, as was also the case in the ATDDG.
The state is denoted by $\textbf{x}=(x_T, y_T, x_A, y_A, x_D, y_D)\in \mathbb{R}^6$; the dynamics are given by \eqref{eq:xT}. 
The co-state is denoted by $\lambda:=(\lambda_{x_A},\lambda_{y_A}, \lambda_{x_D}, \lambda_{y_D}, \lambda_{x_T}, \lambda_{y_T}) \in \mathbb{R}^6$.	 
We consider the operationally relevant cost/payoff function \eqref{eq:costDG2} which is evaluated when $A$ captures $T$ at time $t_f$. 
The terminal time $t_f$ is free and the terminal manifold $\mathcal{S}_c$ is the hyperplane in $\mathbb{R}^6$ 
\begin{align}
\mathcal{S}_c = \left\{ \textbf{x} \ |  \begin{bmatrix}
	1&0&0& 0 & -1 & 0 \\
	0 & 1 & 0 & 0 & 0 & -1 
\end{bmatrix}\begin{bmatrix} x_A \\ y_A \\ x_D \\ y_D \\ x_T \\ y_T  \end{bmatrix} =0_{2\times 1} \right\} 	\label{eq:Hplane}
\end{align}
The Hamiltonian of the differential game is
\begin{align}
  \left.
	\begin{array}{l l}
	\mathcal{H}&=\lambda_{x_A}\cos\chi + \lambda_{y_A}\sin\chi + \lambda_{x_D}\cos\psi \\
	&~~+ \lambda_{y_D}\sin\psi  
	+ \alpha\lambda_{x_T}\cos\phi + \alpha\lambda_{y_T}\sin\phi. 
\end{array}  \right. 
\end{align}
where the Hamiltonian and the dynamics are separable (or decoupled) in the controls $\phi$, $\psi$ and $\chi$. Hence, $\min_{\phi,\psi}  \max_\chi \mathcal{H}= \max_\chi \min_{\phi,\psi}  \mathcal{H}$ and Isaacs' condition holds.

The solution of the game of degree is given in the following Theorem where the synthesis of the state feedback optimal strategies of each player is obtained by employing the ``Two-sided''  or ``Two-person'' PMP \cite{Basar99}, \cite{Lewin94}. Following the construction of the regular solutions the Value function $V$ is explicitly derived and we show that $V$ and $\frac{\partial V}{\partial \textbf{x}}$ are continuous for any $\textbf{x} \in \mathcal{R}_c$; hence, no singular surfaces exist. It is also shown that the Value function globally satisfies the HJI equation in $\mathcal{R}_c$, the winning region of $A$.

We introduce the $A-T$ Apollonius circle 
\begin{align}
\left.
	\begin{array}{l l}
   (x-x_c)^2 + (y-y_c)^2  =r^2
  \end{array}  \right.   \label{eq:ATapp}  
\end{align}
where 
\begin{align}
 \left.
	\begin{array}{l l} 
	x_c=\frac{1}{1-\alpha^2}(x_T-\alpha^2 x_A)  \\
	y_c=\frac{1}{1-\alpha^2}(y_T-\alpha^2 y_A)  \\
	r=\frac{\alpha}{1-\alpha^2}\sqrt{(x_T-x_A)^2+(y_T-y_A)^2}  .
    \end{array}  \right.   \label{eq:CirCenRad}
\end{align}

\begin{theorem}  \label{th:main}
Consider the TAD Capture Differential Game (CDG): \eqref{eq:xT}, \eqref{eq:PCx}-\eqref{eq:costDG3}. Assume that $\textbf{x}\in R_c$. The problem parameter is the speed ratio $0\leq\alpha < 1$. The optimal state feedback strategies of the Target,  the Defender, and the Attacker are given by
\begin{align}
\left.
	 \begin{array}{l l}
  \cos\phi^*= \frac{x^*(\textbf{x})-x_{T}}{\sqrt{(x^*(\textbf{x})-x_{T})^2+(y^*(\textbf{x})-y_{T})^2}} \\
    \sin\phi^*= \frac{y^*(\textbf{x})-y_{T}}{\sqrt{(x^*(\textbf{x})-x_{T})^2+(y^*(\textbf{x})-y_{T})^2}} 
  \end{array}   \right.  \label{eq:OCTar}   \\  
      \left.
	 \begin{array}{l l}
   \cos\psi^*=\frac{x^*(\textbf{x})-x_{D}}{\sqrt{(x^*(\textbf{x})-x_D)^2+(y^*(\textbf{x})-y_D)^2}} \\ 
    \sin\psi^*=\frac{y^*(\textbf{x})-y_D}{\sqrt{(x^*(\textbf{x})-x_D)^2+(y^*(\textbf{x})-y_D)^2}} 
      \end{array}   \right.  \label{eq:OCDef}   \\
       \left.
	 \begin{array}{l l}
  \cos\chi^*=\frac{x^*(\textbf{x})-x_{A}}{\sqrt{(x^*(\textbf{x})-x_A)^2+(y^*(\textbf{x})-y_A)^2}} \\ 
   \sin\chi^*=\frac{y^*(\textbf{x})-y_A}{\sqrt{(x^*(\textbf{x})-x_A)^2+(y^*(\textbf{x})-y_A)^2}}   
     \end{array}   \right.   \label{eq:OC}  
\end{align}  
where the capture point coordinates $(x^*(\textbf{x}),y^*(\textbf{x}))$ is a solution of the system of two equations, \eqref{eq:ATapp}  and the following equation
\begin{align}
\left.
	 \begin{array}{l l}
\  [(x-x_T)^2+(y-y_T)^2] \big[ (x-x_D)(y-y_T) \\
 \ \ \   - (x-x_T)(y-y_D)  - \alpha^2 (x-x_D)(y-y_A) \\  
  \ \ \ +\alpha^2(x-x_A)(y-y_D) \big] ^2  \\
 -   \alpha [(x-x_D)^2+(y-y_D)^2]  \big[(x-x_A)(y-y_T)  \\
 \ \ \ -(x-x_T)(y-y_A)\big]^2 =0.
\end{array}   \right.   \label{eq:Quartic}
\end{align}
The Value function is $C^1$, it satisfies the HJI PDE, and is explicitly given by
\begin{align}
\left.
	 \begin{array}{l l} 
	V(\textbf{x})=& \sqrt{[x^*(\textbf{x})-x_D ]^2+[y^*(\textbf{x}) -y_D ]^2} \\
	&-\frac{1}{\alpha} \sqrt{[x^*(\textbf{x}) -x_T ]^2+[y^*(\textbf{x}) -y_T ]^2}.
\end{array}   \right.  \label{eq:ValueFn}
\end{align}
\end{theorem}

\textit{Proof}. We first determine the optimal control inputs in terms of the co-state variables. This can be directly obtained from Isaacs' Main Equation 1 (ME 1)
\begin{align}
    \min_{\phi,\psi}  \max_{\chi}\mathcal{H} =0 	\label{eq:minmax}
\end{align}
and the optimal control inputs are characterized by the relationships
\begin{align}
\left.
	\begin{array}{l l}
  \cos\chi^*=\frac{\lambda_{x_A}}{\sqrt{\lambda_{x_A}^2+\lambda_{y_A}^2}}, \quad \ \  \sin\chi^*=\frac{\lambda_{y_A}}{\sqrt{\lambda_{x_A}^2+\lambda_{y_A}^2}}  \ \ \ 
  \end{array}  \right.    \label{eq:chico}  \\
  \left.
	\begin{array}{l l}
	\cos\psi^*=-\frac{\lambda_{x_D}}{\sqrt{\lambda_{x_D}^2+\lambda_{y_D}^2}}, \ \  \sin\psi^*=-\frac{\lambda_{y_D}}{\sqrt{\lambda_{x_D}^2+\lambda_{y_D}^2}} 
	\end{array}  \right.  \label{eq:psico}  \\
   \left.
	\begin{array}{l l} 
	 \cos\phi^*=-\frac{\lambda_{x_T}}{\sqrt{\lambda_{x_T}^2+\lambda_{y_T}^2}}, \ \  \sin\phi^*=-\frac{\lambda_{y_T}}{\sqrt{\lambda_{x_T}^2+\lambda_{y_T}^2}}. 
	 \end{array}  \right.   \label{eq:phico} 
\end{align}
 The co-state dynamics are found from the corresponding condition $\dot{\lambda}=-\frac{\partial H}{\partial \textbf{x}}$ which results in: $\dot{\lambda}_{x_A}=\dot{\lambda}_{y_A}=\dot{\lambda}_{x_D}=\dot{\lambda}_{y_D}=\dot{\lambda}_{x_T}=\dot{\lambda}_{y_T}=0$. Therefore, all co-states are constant which means that the optimal headings $\chi^*$, $\psi^*$ and $\phi^*$ are constant as well. The optimal trajectories are straight lines.

Concerning the solution of the attendant Two-Point Boundary Value Problem (TPBVP) in $\mathbb{R}^{12}$ on $0\leq t\leq t_f$, we have six initial states specified by \eqref{eq:xT}. In addition,  six more conditions at the terminal time $t_f$ are needed.
In this respect, define the augmented terminal value function $\Phi_a:\mathbb{R}^6\rightarrow \mathbb{R}^1$ 
\begin{align}
\left.
	\begin{array}{l l}
\Phi_a(\textbf{x}_f)&:= \sqrt{(x_{D_f}-x_{T_f})^2+(y_{D_f}-y_{T_f})^2}  \\
&~~+ \nu_1(x_{A_f}-x_{T_f})+\nu_2(y_{A_f}-y_{T_f})
\end{array}  \right. 
\end{align}
where $\nu_1$ and $\nu_2$ are Lagrange multipliers. The PMP, or Dynamic Programming, directly yields the transversality/terminal co-state conditions
\begin{align}
\lambda(t_f)=\frac{\partial}{\partial \textbf{x}}\Phi_a(\textbf{x}_f)
\end{align}
which results in the following relationships
\begin{align}
&\lambda_{x_A}=\nu_1    \label{eq:Lxd}   \\
&\lambda_{y_A}=\nu_2    \label{eq:Lyd}   \\
&\lambda_{x_D}=\frac{x_{D_f}-x_{T_f}}{\sqrt{(x_{D_f}-x_{T_f})^2+(y_{D_f}-y_{T_f})^2}}    \label{eq:Lxa}   \\
&\lambda_{y_D}=\frac{y_{D_f}-y_{T_f}}{\sqrt{(x_{D_f}-x_{T_f})^2+(y_{D_f}-y_{T_f})^2}}   \label{eq:Lya}   \\
&\lambda_{x_T}=\frac{x_{T_f}-x_{D_f}}{\sqrt{(x_{D_f}-x_{T_f})^2+(y_{D_f}-y_{T_f})^2}} -\nu_1    \label{eq:Lxt}   \\
&\lambda_{y_T}=\frac{y_{T_f}-y_{D_f}}{\sqrt{(x_{D_f}-x_{T_f})^2+(y_{D_f}-y_{T_f})^2}} -\nu_2.    \label{eq:Lyt}   
\end{align}
At this point, we found that \eqref{eq:Lxd}-\eqref{eq:Lyt} together with eqs. \eqref{eq:PCx} yield eight conditions. However, only six conditions are needed and we can eliminate the introduced Lagrange multipliers $\nu_1$ and $\nu_2$ from \eqref{eq:Lxd}-\eqref{eq:Lyd} and  \eqref{eq:Lxt}-\eqref{eq:Lyt} as follows
\begin{align}
&\lambda_{x_A}+\lambda_{x_T}=\frac{x_{T_f}-x_{D_f}}{\sqrt{(x_{D_f}-x_{T_f})^2+(y_{D_f}-y_{T_f})^2}}   \label{eq:Lxad}   \\
&\lambda_{y_A}+\lambda_{y_T}=\frac{y_{T_f}-y_{D_f}}{\sqrt{(x_{D_f}-x_{T_f})^2+(y_{D_f}-y_{T_f})^2}} .   \label{eq:Lyad}   
\end{align}
Then, the six terminal conditions are the two equations in \eqref{eq:PCx} in addition to \eqref{eq:Lxa}-\eqref{eq:Lya} and \eqref{eq:Lxad}-\eqref{eq:Lyad}. Since the terminal time is not fixed, the PMP requirement that the Hamiltonian $\mathcal{H}(\textbf{x}(t),\lambda(t),\chi^*,\psi^*,\phi^*) |_{t_f}\equiv 0$ is used in order to determine $t_f$. Such requirement takes the form of Isaacs' ME 1 
\begin{align}
  \left.
	\begin{array}{r r}
	\lambda_{x_A}\cos\chi^* + \lambda_{y_A}\sin\chi^* + \lambda_{x_D}\cos\psi^*  \\
	 + \lambda_{y_D}\sin\psi^* +\alpha\lambda_{x_T}\cos\phi^* + \alpha\lambda_{y_T}\sin\phi^*\equiv 0.   \label{eq:Hambart}
\end{array}  \right. 
\end{align}
Let $x_T=x_T(t)$, $y_T=y_T(t)$, $x_A=x_A(t)$, $y_A=y_A(t)$, $x_D=x_D(t)$, and $y_D=y_D(t)$ be the instantaneous positions at some time $t<t_f$. From \eqref{eq:xT}, and knowing that the optimal headings of $A$, $T$, and $D$ are constant, we can write
\begin{align}
	&x_{T_f}=x_{T}+\alpha (t_f-t)\cos\phi ^*   \label{eq:xTt}  \\
	&y_{T_f}=y_{T}+\alpha (t_f-t)\sin\phi^*     \label{eq:yTt}  \\
 	&x_{A_f}=x_{A}+(t_f-t)\cos\chi^*     \label{eq:xAt}  \\
	&y_{A_f}=y_A +(t_f-t)\sin\chi^*   \label{eq:yAt} \\
	&x_{D_f}=x_{D}+(t_f-t)\cos\psi^*   \label{eq:xDt}  \\
	&y_{D_f}=y_D+(t_f-t)\sin\psi^*.    	\label{eq:yDt}
\end{align}
According to the terminal condition in \eqref{eq:PCx} we can define 
\begin{align}
x\equiv x_{T_f} = x_{A_f},   \ \ \ \ \ \    y\equiv y_{T_f} = y_{A_f}.
\end{align}
Since $V_T<V_A$ and the optimal headings are constant, then capture of $T$ by $A$ occurs on the $A-T$ Apollonius circle \eqref{eq:ATapp}.

In addition,  the co-state equations   can be written in the following form
\begin{align}
&\lambda_{x_D}=\frac{x_{D_f}-x}{\sqrt{(x_{D_f}-x)^2+(y_{D_f}-y)^2}}  \label{eq:Lxtt} \\
&\lambda_{y_D}=\frac{y_{D_f}-y}{\sqrt{(x_{D_f}-x)^2+(y_{D_f}-y)^2}} \label{eq:Lytt} \\
&\lambda_{x_A}+\lambda_{x_T}=\frac{x-x_{D_f}}{\sqrt{(x_{D_f}-x)^2+(y_{D_f}-y)^2}}  \label{eq:Lxadt}   \\
&\lambda_{y_A}+\lambda_{y_T}=\frac{y-y_{D_f}}{\sqrt{(x_{D_f}-x)^2+(y_{D_f}-y)^2}}  .  \label{eq:Lyadt}   
\end{align}
Substituting the co-states $\lambda_{x_D}$ and $\lambda_{y_D}$  into \eqref{eq:psico} we obtain the optimal Defender heading in terms of $(x_{D_f},y_{D_f})$ which is the terminal state of the Defender
\begin{align}
\left.
	 \begin{array}{l l} 
\cos\psi^*=\frac{x-x_{D_f}}{\sqrt{(x-x_{D_f})^2+(y-y_{D_f})^2}}   \\  \sin\psi^*=\frac{y-y_{D_f}}{\sqrt{(x-x_{D_f})^2+(y-y_{D_f})^2}}. 
\end{array}   \right.   \label{eq:phigeo} 
\end{align}
Here we recall again that the optimal headings are constant, therefore \eqref{eq:phigeo} can be equivalently written in terms of the instantaneous state $(x_D,y_D)$ as it is shown next
\begin{align}
\left.
	 \begin{array}{l l} 
 \cos\psi^*=\frac{x-x_{D}}{\sqrt{(x-x_{D})^2+(y-y_{D})^2}} \\
    \sin\psi^*=\frac{y-y_{D}}{\sqrt{(x-x_{D})^2+(y-y_{D})^2}}. 
\end{array}   \right.  \label{eq:phia}
\end{align}
Additionally, since the optimal headings of $A$ and $T$ are constant and $A$ captures $T$ at the interception point $I:(x,y)$, then we can write the optimal headings of the Attacker and the Target in terms of their corresponding instantaneous states $(x_A,y_A)$ and $(x_T,y_T)$, respectively
\begin{align}
\left.
	 \begin{array}{l l} 
  \cos\chi^*=\frac{x-x_{A}}{\sqrt{(x-x_{A})^2+(y-y_{A})^2}}  \\ 
   \sin\chi^*=\frac{y-y_A}{\sqrt{(x-x_{A})^2+(y-y_{A})^2}}  
  \end{array}   \right.  \label{eq:chi}   
\end{align}
 and
 \begin{align}
  \left.
	 \begin{array}{l l} 
	  \cos\phi^*=\frac{x-x_T}{\sqrt{(x-x_T)^2+(y-y_T)^2}} \\  
	   \sin\phi^*=\frac{y-y_T}{\sqrt{(x-x_T)^2+(y-y_T)^2}}. 
	\end{array}   \right.  \label{eq:psi} 
\end{align}
%
From  \eqref{eq:chico} and \eqref{eq:chi} we have that the Attacker's states and co-states satisfy
\begin{align}
\left.
	 \begin{array}{l l} 
  \frac{\lambda_{x_A}}{\sqrt{\lambda_{x_A}^2+\lambda_{y_A}^2}}=\frac{x-x_{A}}{\sqrt{(x-x_A)^2+(y-y_A)^2}} \\
   \frac{\lambda_{y_A}}{\sqrt{\lambda_{x_A}^2+\lambda_{y_A}^2}}=\frac{y-y_A}{\sqrt{(x-x_A)^2+(y-y_A)^2}}.   
   \end{array}   \right.   \label{eq:Acostate} 
\end{align}
Similarly, using  \eqref{eq:phico} and \eqref{eq:psi} we can obtain the corresponding relationships for the Target
\begin{align}
\left.
	 \begin{array}{l l} 
  -\frac{\lambda_{x_T}}{\sqrt{\lambda_{x_T}^2+\lambda_{y_T}^2}}=\frac{x-x_T}{\sqrt{(x-x_T)^2+(y-y_T)^2}} \\ 
  -\frac{\lambda_{y_T}}{\sqrt{\lambda_{x_T}^2+\lambda_{y_T}^2}}=\frac{y-y_T}{\sqrt{(x-x_T)^2+(y-y_T)^2}}. 
  \end{array}   \right.   \label{eq:Dcostate}  
\end{align}
We have four equations \eqref{eq:Lxadt}, \eqref{eq:Lyadt}, \eqref{eq:Acostate}, and \eqref{eq:Dcostate} in the four unknowns $\lambda_{x_A}$, $\lambda_{y_A}$, $\lambda_{x_T}$, and $\lambda_{y_T}$. Making the corresponding substitutions between these equations we can write the following
\begin{align}
\left.
	 \begin{array}{l l}
  \sqrt{\lambda_{x_A}^2+\lambda_{y_A}^2} &= \frac{\sqrt{(x-x_A)^2+(y-y_A)^2}}{\sqrt{(x_{D_f}-x)^2+(y_{D_f}-y)^2}} \\
  &~~ \times \frac{x-x_{D_f} - \frac{x-x_T}{y-y_T}(y-y_{D_f})}{x-x_A-\frac{x-x_T}{y-y_T}(y-y_A)}    \\
  \sqrt{\lambda_{x_T}^2+\lambda_{y_T}^2} &= \frac{\sqrt{(x-x_T)^2+(y-y_T)^2}}{\sqrt{(x_{D_f}-x)^2+(y_{D_f}-y)^2}} \\
  &~~ \times \frac{x-x_{D_f} - \frac{x-x_A}{y-y_A}(y-y_{D_f})}{\frac{x-x_A}{y-y_A}(y-y_T)-(x-x_T)}      
\end{array}   \right.   \nonumber
\end{align}
Substituting the previous equations back into \eqref{eq:Acostate} and \eqref{eq:Dcostate} we obtain the solution
\begin{align}
\left.
	 \begin{array}{l l}
  \lambda_{x_A}= \frac{x-x_A}{\sqrt{(x_{D_f}-x)^2+(y_{D_f}-y)^2}} \cdot \frac{x-x_{D_f} - \frac{x-x_T}{y-y_T}(y-y_{D_f})}{x-x_A-\frac{x-x_T}{y-y_T}(y-y_A)}     \\  
  \lambda_{y_A}=  \frac{y-y_A}{\sqrt{(x_{D_f}-x)^2+(y_{D_f}-y)^2}} \cdot \frac{x-x_{D_f} - \frac{x-x_T}{y-y_T}(y-y_{D_f})}{x-x_A-\frac{x-x_T}{y-y_T}(y-y_A)}    \\    
  \lambda_{x_T}= \frac{x_T-x}{\sqrt{(x_{D_f}-x)^2+(y_{D_f}-y)^2}} \cdot \frac{x-x_{D_f} - \frac{x-x_A}{y-y_A}(y-y_{D_f})}{\frac{x-x_A}{y-y_A}(y-y_T) - (x-x_T)}  \\  
  \lambda_{y_T}= \frac{y_T-y}{\sqrt{(x_{D_f}-x)^2+(y_{D_f}-y)^2}} \cdot \frac{x-x_{D_f} - \frac{x-x_A}{y-y_A}(y-y_{D_f})}{\frac{x-x_A}{y-y_A}(y-y_T) - (x-x_T)}. 
\end{array}   \right.   \label{eq:Lydt}
\end{align}
These equations together with  \eqref{eq:Lxtt} and \eqref{eq:Lytt} specify the co-states in terms of the instantaneous state components $(x_A,y_A,x_T,y_T)$ but also in terms of the Defender's terminal position $(x_{D_f},y_{D_f})$. At this point we recall that  the dynamical equations were normalized by $V_A$; hence, the terminal time satisfies 
\begin{align}
\left.
	 \begin{array}{l l} 
	t_f-t&=\overline{AI}= \frac{1}{\alpha} \overline{TI} \\
	&= \frac{1}{\alpha} \sqrt{(x-x_T)^2+(y-y_T)^2}
\end{array}   \right.  \label{eq:bart} 
\end{align}
which can be used to eliminate the dependence of equations \eqref{eq:Lxtt}-\eqref{eq:Lytt} and \eqref{eq:Lydt} on the Defender's terminal position. The goal is to determine expressions for each co-state that do not depend on terminal states but on the instantaneous states in order to obtain a state feedback strategy.  This can be achieved by using \eqref{eq:xDt}, \eqref{eq:yDt}, \eqref{eq:phia}, and \eqref{eq:bart}  in order to obtain the following equations    
\begin{align}
	\left.
 \begin{array}{l l}
     x_{D_f} \!\!
      &= x_D  \\
      &~~  + \frac{1}{\alpha} \sqrt{(x-x_T)^2+(y-y_T)^2}  \frac{x-x_D}{\sqrt{(x-x_D)^2+(y-y_D)^2}}   \\
	   y_{D_f} \!\! 
      &= y_D \\
      &~~ + \frac{1}{\alpha} \sqrt{(x-x_T)^2+(y-y_T)^2}  \frac{y-y_D}{\sqrt{(x-x_D)^2+(y-y_D)^2}} .  	\end{array}   \right.  \label{eq:yTopta}  
\end{align}
Now we can use  \eqref{eq:yTopta} in equations \eqref{eq:Lxtt}-\eqref{eq:Lytt} and \eqref{eq:Lydt} and obtain the expressions of the co-states in terms of the state $\textbf{x}$ and the point of interception coordinates $(x,y)$:
\begin{align}
\left.
	 \begin{array}{l l}
     \lambda_{x_D}= \frac{x_D-x}{\sqrt{(x_D-x)^2+(y_D-y)^2}}  \\
	 \lambda_{y_D}= \frac{y_D-y}{\sqrt{(x_D-x)^2+(y_D-y)^2}}  \\
	 \lambda_{x_A}= \frac{x-x_A}{\sqrt{(x_D-x)^2+(y_D-y)^2}} \cdot \frac{x-x_D - \frac{x-x_T}{y-y_T}(y-y_D)}{x-x_A-\frac{x-x_T}{y-y_T}(y-y_A)}   \\  
  \lambda_{y_A}=  \frac{y-y_A}{\sqrt{(x_D-x)^2+(y_D-y)^2}} \cdot \frac{x-x_D - \frac{x-x_T}{y-y_T}(y-y_D)}{x-x_A-\frac{x-x_T}{y-y_T}(y-y_A)}    \\    
  \lambda_{x_T}= \frac{x_T-x}{\sqrt{(x_D-x)^2+(y_D-y)^2}} \cdot \frac{x-x_D - \frac{x-x_A}{y-y_A}(y-y_D)}{\frac{x-x_A}{y-y_A}(y-y_T) - (x-x_T)}  \\  
  \lambda_{y_T}= \frac{y_T-y}{\sqrt{(x_D-x)^2+(y_D-y)^2}} \cdot \frac{x-x_D - \frac{x-x_A}{y-y_A}(y-y_D)}{\frac{x-x_A}{y-y_A}(y-y_T) - (x-x_T)}  .   \\   
	\end{array}   \right.   \label{eq:costateopta}  
\end{align}
We now turn our attention to equation \eqref{eq:Hambart}, where we use \eqref{eq:phia}-\eqref{eq:psi} for the optimal controls and  \eqref{eq:costateopta} for the co-states and we have that
\begin{align}
\left.
	 \begin{array}{l l}
& \frac{(x-x_A)^2+(y-y_A)^2}{\sqrt{(x-x_A)^2+(y-y_A)^2}} \cdot \frac{x-x_D - \frac{x-x_T}{y-y_T}(y-y_D)}{x-x_A-\frac{x-x_T}{y-y_T}(y-y_A)} \\
& - \frac{(x_D-x)^2+(y_D-y)^2}{\sqrt{(x_D-x)^2+(y_D-y)^2}}   \\
& - \alpha \frac{(x-x_T)^2+(y-y_T)^2}{\sqrt{(x-x_T)^2+(y-y_T)^2}}  \cdot \frac{x-x_D - \frac{x-x_A}{y-y_A}(y-y_D)}{\frac{x-x_A}{y-y_A}(y-y_T) - (x-x_T)}   =0   
\end{array}   \right.   \nonumber
\end{align}
 after the common term $\frac{1}{\sqrt{(x_D-x)^2+(y_D-y)^2}}$ has been canceled. 
Using the relationships in \eqref{eq:bart}  we can write
\begin{align}
\left.
	 \begin{array}{l l}
 \  \frac{1}{\alpha}\sqrt{(x-x_T)^2+(y-y_T)^2}  \frac{(x-x_D)(y-y_T)-(x-x_T)(y-y_D)}{(x-x_A)(y-y_T)-(x-x_T)(y-y_A)}  \\
 \ - \alpha \sqrt{(x-x_T)^2+(y-y_T)^2}  \frac{(x-x_D)(y-y_A)-(x-x_A)(y-y_D)}{(x-x_A)(y-y_T)-(x-x_T)(y-y_A)}   \\
 =   \sqrt{(x-x_D)^2+(y-y_D)^2}  . 
\end{array}   \right.   \nonumber
\end{align}
Grouping terms and multiplying both sides of the previous equation by $\alpha [(x-x_A)(y-y_T) -(x-x_T)(y-y_A)]$ we obtain
\begin{align}
\left.
	 \begin{array}{l l}
\  \sqrt{(x-x_T)^2+(y-y_T)^2} \big[ (x-x_D)(y-y_T) \\
 \ \ \   - (x-x_T)(y-y_D)  - \alpha^2 (x-x_D)(y-y_A) \\  
  \ \ \ +\alpha^2(x-x_A)(y-y_D) \big]   \\
 =   \alpha \sqrt{(x-x_D)^2+(y-y_D)^2}  \big[(x-x_A)(y-y_T)  \\
 \ \ \ -(x-x_T)(y-y_A)\big]
\end{array}   \right.   \label{eq:xyoptimal}
\end{align}
which can be written as in \eqref{eq:Quartic}. The optimal interception point $I^*:(x^*,y^*)$ is obtained by solving the systems of polynomial equations \eqref{eq:ATapp}  and \eqref{eq:Quartic}. Therefore, the state feedback optimal strategies have been obtained and the value function is determined; it is $C^1$ and is explicitly given by \eqref{eq:ValueFn}. The gradient of the Value function, which is shown in \eqref{eq:pVpx} below, is well defined in the region $\mathcal{R}_c$.
In order to obtain $\frac{\partial V(\textbf{x})}{\partial\textbf{x}}$ we first compute the following term 
from equation \eqref{eq:ATapp}
  \begin{align}
\left.
	 \begin{array}{l l} 
	\frac{dy^*}{dx^*} = -\frac{x^*-x_c}{y^*-y_c}
\end{array}   \right.  \label{eq:dervyx}
\end{align}
 where $x_c,y_c,r$ are given by \eqref{eq:CirCenRad}.

The gradient of $V(\textbf{x})$ is $\frac{\partial V}{\partial \textbf{x}}=[\frac{\partial V}{\partial x_i}\! +\! \frac{dV}{dx^*}\cdot\frac{dx^*}{dx_i}, \ \ \ \frac{\partial V}{\partial y_i}\! +\! \frac{dV}{dx^*}\cdot\frac{dx^*}{dy_i}]^T$ for $i=A,T,D$.
We start by determining the following term 
\begin{align}
\left.
	 \begin{array}{l l} 
\frac{d V}{dx^*} =&   \frac{x^* - x_D - (y^*-y_D)\frac{x^*-x_c}{y^*-y_c}}{\sqrt{(x^*-x_D)^2+(y^*-y_D)^2}}   \\
 &-  \frac{x^* - x_T - (y^*-y_T)\frac{x^*-x_c}{y^*-y_c}}{\alpha \sqrt{(x^*-x_T)^2+(y^*-y_T)^2}}  \label{eq:Derxop}
\end{array}   \right. 
\end{align}
where $y^*=y_c + \sqrt{r^2-(x^*-x_c)^2}$.
 
 Let us now write \eqref{eq:xyoptimal} as follows
\begin{align}
\left.
	 \begin{array}{l l} 
\   \frac{ (x-x_D)(y-y_T) - (x-x_T)(y-y_D)  - \alpha^2 [(x-x_D)(y-y_A) -(x-x_A)(y-y_D)] }{\sqrt{(x-x_D)^2+(y-y_D)^2} }   \\
 =  \frac{ \alpha^2[(x-x_A)(y-y_T)  -(x-x_T)(y-y_A)] }{\alpha \sqrt{(x-x_T)^2+(y-y_T)^2}}
\end{array}   \right.  \label{eq:QxyP}
\end{align}
where the right hand side of this equation was multiplied and divided by $\alpha$. Let us add and subtract the term $(x-x_T)(y-y_T)$ to the numerator of the right hand side of the equation and, by rearranging terms, we have 
\begin{align}
\left.
	 \begin{array}{l l} 
 \  \frac{ (x-x_D)(y-y_T-\alpha^2 y +\alpha^2 y_A) - (y-y_D)(x-x_T-\alpha^2 x +\alpha^2 x_A)  ] }{\sqrt{(x-x_D)^2+(y-y_D)^2} }   \\
 -  \frac{ (x-x_T)(y-y_T-\alpha^2 y +\alpha^2 y_A) - (y-y_T)(x-x_T-\alpha^2 x +\alpha^2 x_A)  }{\alpha \sqrt{(x-x_T)^2+(y-y_T)^2}} = 0.
\end{array}   \right.  \label{eq:QxyP2}
\end{align}
Dividing both sides of \eqref{eq:QxyP2} by $1-\alpha^2$ we obtain
\begin{align}
\left.
	 \begin{array}{l l} 
 \  \frac{ (x-x_D) \frac{(1-\alpha^2) y-(y_T -\alpha^2 y_A)}{1-\alpha^2} - (y-y_D)\frac{(1-\alpha^2) x-(x_T -\alpha^2 x_A)}{1-\alpha^2}  }{\sqrt{(x-x_D)^2+(y-y_D)^2} }   \\
 -  \frac{ (x-x_T)\frac{(1-\alpha^2) y-(y_T -\alpha^2 y_A)}{1-\alpha^2} - (y-y_T)\frac{(1-\alpha^2) x-(x_T -\alpha^2 x_A)}{1-\alpha^2}  }{\alpha \sqrt{(x-x_T)^2+(y-y_T)^2}} = 0.
\end{array}   \right. \nonumber 
\end{align}
Using the definitions in \eqref{eq:CirCenRad} we have that
\begin{align}
\left.
	 \begin{array}{l l} 
 \  \frac{ (x-x_D) (y-y_c)- (y-y_D)(x-x_c) }{\sqrt{(x-x_D)^2+(y-y_D)^2} }   \\
 -  \frac{ (x-x_T)(y-y_c) - (y-y_T)(x-x_c) }{\alpha \sqrt{(x-x_T)^2+(y-y_T)^2}} = 0.
\end{array}   \right.  \label{eq:QxyP4}
\end{align}
Finally, dividing both sides  of \eqref{eq:QxyP4} by $y-y_c$ we obtain
\begin{align}
\left.
	 \begin{array}{l l} 
 \  \frac{ (x-x_D)- (y-y_D)\frac{x-x_c}{y-y_c} }{\sqrt{(x-x_D)^2+(y-y_D)^2} }  
 -  \frac{ (x-x_T) - (y-y_T)\frac{x-x_c}{y-y_c} }{\alpha \sqrt{(x-x_T)^2+(y-y_T)^2}} = 0 .
\end{array}   \right.  \label{eq:QxyP5}
\end{align}
We have found that $\frac{d V}{dx^*}$ is equal to the left hand side of  \eqref{eq:QxyP5} where $y^*=y_c + \sqrt{r^2-(x^*-x_c)^2}$ and $x^*$ is a solution of \eqref{eq:Quartic}, equivalently, of \eqref{eq:QxyP5}. Therefore, we have that 
 \begin{align}
\left.
	 \begin{array}{l l} 
\frac{d V}{dx^*} =  0.
\end{array}   \right.   \nonumber
\end{align}
We now substitute $y^*=y_c + \sqrt{r^2-(x^*-x_c)^2}$ in \eqref{eq:ValueFn}
\begin{align}
\left.
	 \begin{array}{l l} 
	V(\textbf{x})= \\
	 \sqrt{[x^*(\textbf{x})-x_D ]^2+[y_c + \sqrt{r^2-(x^*(\textbf{x})-x_c)^2}-y_D ]^2} \\
	-\frac{1}{\alpha} \sqrt{[x^*(\textbf{x}) -x_T ]^2+[y_c + \sqrt{r^2-(x^*(\textbf{x})-x_c)^2} -y_T ]^2}.
\end{array}   \right.  \nonumber
\end{align}

Since we have shown that $\frac{d V}{dx^*} =  0$, then, the gradient of $V(\textbf{x})$ is simplified and given by $\frac{\partial V}{\partial \textbf{x}}=[\frac{\partial V}{\partial x_A} \ \  \frac{\partial V}{\partial y_A} \ \ \frac{\partial V}{\partial x_T} \ \  \frac{\partial V}{\partial y_T} \ \ \frac{\partial V}{\partial x_D} \ \  \frac{\partial V}{\partial y_D} ]^T$. We now obtain the following
 \begin{align}
\left.
	 \begin{array}{l l} 
	\frac{\partial x_c}{\partial x_T} =\frac{\partial y_c}{\partial y_T} = \frac{1}{1-\alpha^2}  \\
	\frac{\partial x_c}{\partial x_A} =\frac{\partial y_c}{\partial y_A} = -\frac{\alpha^2}{1-\alpha^2}  \\
	\frac{\partial r^2}{\partial x_T} =  \frac{2\alpha^2}{(1-\alpha^2)^2} (x_T-x_A)    \\
	\frac{\partial r^2}{\partial y_T} =  \frac{2\alpha^2}{(1-\alpha^2)^2} (y_T-y_A)    \\
	\frac{\partial r^2}{\partial x_A} =  -\frac{2\alpha^2}{(1-\alpha^2)^2} (x_T-x_A)    \\
	\frac{\partial r^2}{\partial y_A} = - \frac{2\alpha^2}{(1-\alpha^2)^2} (y_T-y_A).
\end{array}   \right.  \label{eq:pdxcycr}
\end{align}
The previous expressions are useful to compute the gradient of $V(\textbf{x})$. Then, we can write
 \begin{align}
\left.
	 \begin{array}{l l} 
	\frac{\partial V}{\partial x_A} =  \\
	 - \frac{\alpha^2 [y_c + \sqrt{r^2-(x^*-x_c)^2} -y_D] [ x_T-x_A + (1-\alpha^2) (x^*-x_c)] }{(1-\alpha^2)^2 \sqrt{r^2-(x^*-x_c)^2}  \sqrt{(x^*-x_D )^2+(y_c + \sqrt{r^2-(x^*-x_c)^2}-y_D )^2}}  \\
	  + \frac{\alpha [y_c + \sqrt{r^2-(x^*-x_c)^2} -y_T] [ x_T-x_A + (1-\alpha^2) (x^*-x_c)] }{(1-\alpha^2)^2 \sqrt{r^2-(x^*-x_c)^2}  \sqrt{(x^*-x_T)^2+(y_c + \sqrt{r^2-(x^*-x_c)^2}-y_T )^2}}
\end{array}   \right.  \label{eq:pvxA}
\end{align}
which can be simplified in the following form
 \begin{align}
\left.
	 \begin{array}{l l} 
	\frac{\partial V}{\partial x_A} &=  - \frac{\alpha^2 [y^*-y_D] [ x_T-x_A + (1-\alpha^2) (x^*-x_c)] }{(1-\alpha^2)^2 (y^*-y_c)  \sqrt{(x^*-x_D)^2+(y^*-y_D)^2}}  \\
	 & ~~ + \frac{\alpha [y^* -y_T] [ x_T-x_A + (1-\alpha^2) (x^*-x_c)] }{(1-\alpha^2)^2  (y^*-y_c)  \sqrt{(x^*-x_T)^2+(y^*-y_T)^2}}  \\
	 &=  - \frac{\alpha^2 (y^*-y_D) (x^* -x_A) }{(1-\alpha^2) (y^*-y_c)  \sqrt{(x^*-x_D)^2+(y^*-y_D )^2}}  \\
	  & ~~ + \frac{\alpha (y^*-y_T) (x^* -x_A) }{(1-\alpha^2) (y^*-y_c)  \sqrt{(x^*-x_T)^2+(y^*-y_T )^2}}  \\
	  & = \alpha \frac{x^* -x_A}{(1-\alpha^2) (y^*-y_c)}  \big( \frac{y^*-y_T}{\overline{TI}} - \alpha \frac{y^*-y_D}{\overline{DI}} \big)
\end{array}   \right.  \label{eq:pvxA2}
\end{align}
where $\overline{TI} = \sqrt{(x^*-x_T)^2+(y^*-y_T )^2}$ and $\overline{DI} = \sqrt{(x^*-x_D)^2+(y^*-y_D)^2}$. 
We can now write the partial derivatives of the Value function with respect to each component of the state $\textbf{x}$ and they are given by
\begin{align}
\left.
	 \begin{array}{l l} 
	 \frac{\partial V}{\partial x_A} &= \alpha \frac{x^* -x_A}{(1-\alpha^2) (y^*-y_c)}  \big( \frac{y^*-y_T}{\overline{TI}} - \alpha \frac{y^*-y_D}{\overline{DI}} \big) \\
	  \frac{\partial V}{\partial y_A} &= \alpha \frac{y^* -y_A}{(1-\alpha^2) (y^*-y_c)}  \big( \frac{y^*-y_T}{\overline{TI}} - \alpha \frac{y^*-y_D}{\overline{DI}} \big)  \\  
	   \frac{\partial V}{\partial x_T} &=  \frac{x^* -x_T}{(1-\alpha^2) (y^*-y_c)}  \big( \frac{y^*-y_D}{\overline{DI}} -  \frac{y^*-y_A}{\overline{AI}} \big)   \\  
	    \frac{\partial V}{\partial y_T} &=  \frac{y^* -y_T}{(1-\alpha^2) (y^*-y_c)}  \big( \frac{y^*-y_D}{\overline{DI}} -  \frac{y^*-y_A}{\overline{AI}} \big)  \\  
	     \frac{\partial V}{\partial x_D} &= - \frac{x^*-x_D}{\overline{DI}} \\
	      \frac{\partial V}{\partial x_D} &= - \frac{y^*-y_D}{\overline{DI}} 
\end{array}   \right.  \label{eq:pVpx} 
\end{align}
where $\overline{AI} = \sqrt{(x^*-x_A)^2+(y^*-y_A )^2}$. 

Finally, we show that the $V(\textbf{x})$ is the solution of the HJI equation
$-\frac{\partial V}{\partial t} =\frac{\partial V}{\partial \textbf{x}}\cdot  \textbf{f}(\textbf{x},\chi^*,\psi^*,\phi^*) + g(t,\textbf{x},\chi^*,\psi^*,\phi^*) $.
Note that in this problem $\frac{\partial V}{\partial t}=0$ and  $g(t,\textbf{x},\chi^*,\psi^*,\phi^*)=0$. 
The HJI equation for the CDG is then given by
\begin{align}
\left.
	 \begin{array}{l l} 
	&\frac{\partial V}{\partial \textbf{x}} \cdot \textbf{f}(\textbf{x},\chi^*,\psi^*,\phi^*)  \\
& ~~ =  -\frac{(x^*-x_D)^2}{\overline{DI}} -  \frac{(y^*-y_D)^2}{\overline{DI}}   \\
& ~~~~ +  \alpha \frac{(x^* -x_A)^2}{(1-\alpha^2) (y^*-y_c) \overline{AI}}  \big( \frac{y^*-y_T}{\overline{TI}} - \alpha \frac{y^*-y_D}{\overline{DI}} \big)  \\
& ~~~~ +  \alpha \frac{(y^* -y_A)^2}{(1-\alpha^2) (y^*-y_c) \overline{AI}}  \big( \frac{y^*-y_T}{\overline{TI}} - \alpha \frac{y^*-y_D}{\overline{DI}} \big)  \\
& ~~~~ + \alpha \frac{(x^* -x_T)^2}{(1-\alpha^2) (y^*-y_c)\overline{TI}}  \big( \frac{y^*-y_D}{\overline{DI}} -  \frac{y^*-y_A}{\overline{AI}} \big)  \\
& ~~~~ + \alpha \frac{(y^* -y_T)^2}{(1-\alpha^2) (y^*-y_c)\overline{TI}}  \big( \frac{y^*-y_D}{\overline{DI}} -  \frac{y^*-y_A}{\overline{AI}} \big) 
\end{array}   \right. 
\end{align}
Since $A$ and $D$ have the same speed, then the distances satisfy $\overline{AI}=\overline{DI}$. Similarly, $\overline{TI}=\alpha \overline{AI}$ and we have that
\begin{align}
\left.
	 \begin{array}{l l} 
	&\frac{\partial V}{\partial \textbf{x}} \cdot \textbf{f}(\textbf{x},\chi^*,\psi^*,\phi^*)  \\
& ~~ =  - 1  +  \frac{\alpha \overline{AI}}{(1-\alpha^2) (y^*-y_c) }  \big( \frac{y^*-y_T}{\overline{TI}} - \alpha \frac{y^*-y_D}{\overline{DI}} \big)  \\
& ~~~~  + \frac{\alpha \overline{TI}}{(1-\alpha^2) (y^*-y_c) }  \big( \frac{y^*-y_D}{\overline{DI}} -  \frac{y^*-y_A}{\overline{AI}} \big)  \\
& ~~ = -1 +  \frac{\alpha }{(1-\alpha^2) (y^*-y_c) }  \big[ \frac{1}{\alpha}(y^*-y_T)  \\
& ~~~~~~ - \alpha (y^*-y_D) + \alpha (y^*-y_D) -  \alpha (y^*-y_A) \big]  \\
& ~~ = - 1+ \frac{1 }{(1-\alpha^2) (y^*-y_c) }  \big[ y^*-y_T -\alpha^2(y^*-y_A) \big]  \\
& ~~ = - 1+ \frac{1 }{(1-\alpha^2) (y^*-y_c) }  \big[(1-\alpha^2) y^*  - (1-\alpha^2) y_c \big]   \\
& ~~ = - 1+ \frac{(1-\alpha^2) (y^*-y_c) }{(1-\alpha^2) (y^*-y_c) }   \\
& ~~ = 0.
\end{array}   \right. 
\end{align}
In summary, when $\textbf{x}\in\mathcal{R}_c$, state feedback optimal strategies of the three agents were synthesized and the Value function was obtained. 
It was also shown that the Value function is $C^1$ and that it is the solution of the HJI equation.   \  $\square$

In order to obtain the optimal capture coordinates $(x^*(\textbf{x}),y^*(\textbf{x}))$ the system of equations \eqref{eq:ATapp}  and \eqref{eq:Quartic} needs to be solved. The following result provides an equivalent and more compact expression to determine the optimal coordinates.

\begin{corollary}
The optimal capture coordinates are given by $x^*(\textbf{x})=x_c+r\cos\omega^*$ and $y^*(\textbf{x})=y_c+r\sin\omega^*$ where $\omega^*$ is a solution of the following polynomial equation
\begin{align}
   \left.
	 \begin{array}{l l}
   \frac{r(a_0 - ib_0)}{4}v^6 + \frac{a_1 - ib_1}{4}v^5 + r(a_2 - ib_2)v^4 + \frac{r^2a_3}{2}v^3  \\
   + r(a_2 + ib_2)v^2 +\frac{a_1 + ib_1}{4}v + \frac{r(a_0 + ib_0)}{4} =0
	\end{array}  \label{eq:complex}  \right.
\end{align}
where $v=e^{iw}$ and 
\begin{align}
   \left.
	 \begin{array}{l l}
a_0 &= [(y_c\!-\!y_D)^2-(x_c\!-\!x_D)^2](x_c-x_A)\\
&~~ +[(x_c\!-\!x_A)^2-(y_c\!-\!y_A)^2](x_c-x_D)  \\
&~~ +2(y_c-y_D)(y_c-y_A)(x_A-x_D)  \\
b_0&=[(y_c\!-\!y_D)^2-(x_c\!-\!x_D)^2](y_c-y_A)\\
&~~ +[(x_c\!-\!x_A)^2-(y_c\!-\!y_A)^2](y_c-y_D)  \\
&~~ -2(x_c-x_D)(x_c-x_A)(y_A-y_D)  \\
a_1 &= 2[ (x_c\!-\!x_A)^2(y_c\!-\!y_D)^2- (x_c\!-\!x_D)^2(y_c\!-\!y_A)^2] \\
&~~  +r^2[(x_c\!-\!x_A)^2 \!-\!(x_c\!-\!x_D)^2 \!+\! (y_c\!-\!y_D)^2\!-\!(y_c\!-\!y_A)^2]  \\
b_1 &= 2[(x_c\!-\!x_D)^2 +(y_c\!-\!y_D)^2+r^2](x_c\!-\!x_A)(y_c\!-\!y_A)  \\
&~~  -2[(x_c\!-\!x_A)^2 +(y_c\!-\!y_A)^2+r^2](x_c\!-\!x_D)(y_c\!-\!y_D)  \\
a_2 &= \frac{[(x_c\!-\!x_D)^2+3(y_c\!-\!y_D)^2](x_c-x_A)-[(x_c\!-\!x_A)^2+3(y_c\!-\!y_A)^2](x_c-x_D)}{4} \\
 &~~ + \frac{(y_c-y_D)(y_c-y_A)(x_D-x_A)}{2} \\
 b_2 &= \frac{[3(x_c\!-\!x_D)^2+(y_c\!-\!y_D)^2](y_c-y_A)-[3(x_c\!-\!x_A)^2+(y_c\!-\!y_A)^2](y_c-y_D)}{4} \\
 &~~ + \frac{(x_c-x_D)(x_c-x_A)(y_D-y_A)}{2} \\
 a_3 & = (x_c\!-\!x_D)^2+(y_c\!-\!y_D)^2-(x_c\!-\!x_A)^2-(y_c\!-\!y_A)^2.
 \end{array}  \label{eq:CoeffPoly}  \right.
\end{align}
\end{corollary}

\section{Game of Degree on the Barrier Surface}  \label{sec:Barrier}
In this section we show that both Games of Degree in the TAD differential game: the CDG and the ATDDG, provide the same solution when  $\textbf{x}\in\mathcal{B}$.

\begin{theorem}
On the Barrier surface, that is, when $\textbf{x}\in\mathcal{B}$, the CDG and the ATDDG provide the same solution and the same Value of the game.
\end{theorem}
\textit{Proof}. When the state of the system is such that $\textbf{x}\in\mathcal{B}$, then it holds that the $A-T$ Apollonius circle is tangent to the orthogonal bisector of the segment $\overline{AD}$. Let the tangent point be $I:(x,y)$. It is easy to show that the line $\overline{CI}$ is parallel to the line $\overline{AD}$, where the point $C$ is the center of the $A-T$ Apollonius circle. Equivalently, one can write 
\begin{align}
\left.
	 \begin{array}{l l} 
	\frac{x-x_c}{y-y_c} = \frac{x_D-x_A}{y_D-y_A}
\end{array}   \right.   \label{eq:slopes}
\end{align}
Then, we are going to show that $\omega^*=\rho$ in this case, where $\rho$ is the LOS angle form $A$ to $D$; this is equivalent to selecting the interception point $I:(x,y)$, the tangent point. It will be shown that this aimpoint choice solves the standing equation \eqref{eq:Quartic}. Note that the tangent point also solves \eqref{eq:ATapp} since it is a point on the Apollonius circle.
 We will use \eqref{eq:Quartic} in the form of \eqref{eq:QxyP5} and, since equation \eqref{eq:slopes} holds for this selection of interception point we can write the left hand side of \eqref{eq:QxyP5} as follows
 \begin{align}
\left.
	 \begin{array}{l l} 
 \  \frac{ (x-x_D)- (y-y_D)\frac{x_D-x_A}{y_D-y_A} }{\sqrt{(x-x_D)^2+(y-y_D)^2} }  
 -  \frac{ (x-x_T) - (y-y_T)\frac{x_D-x_A}{y_D-y_A} }{\alpha \sqrt{(x-x_T)^2+(y-y_T)^2}} 
\end{array}   \right.  \label{eq:xinB}
\end{align}
we note that, since point $I$ is located both on the orthogonal bisector of the segment $\overline{AD}$ and on the $A-T$ Apollonius circle, then the terminal position of both $T$ and $D$ is equal to $I$. Hence, their traveled distance satisfy $\overline{TI}=\alpha\overline{DI}$ and \eqref{eq:xinB} can be written as
\begin{align}
\left.
	 \begin{array}{l l} 
& \frac{1}{\alpha^2 (y_D-y_A) \sqrt{(x-x_D)^2+(y-y_D)^2} }  \\
& \ \times \big( \alpha^2[(x-x_D)(y_D-y_A)- (y-y_D)(x_D-x_A) ] \\
 &~~  - (x-x_T)(y_D-y_A) + (y-y_T)(x_D-x_A) \big).
\end{array}   \right.  \label{eq:xinB2}
\end{align}
Grouping like terms, the expression in  \eqref{eq:xinB2} can be written as follows
\begin{align}
\left.
	 \begin{array}{l l} 
 \frac{-1}{\alpha^2 (y_D-y_A) \sqrt{(x-x_D)^2+(y-y_D)^2} }  \\
 \ \times \big( (1\!-\!\alpha^2)[x(y_D -y_A) +y(x_A-x_D)] + y_T(x_D-x_A) \\
 ~~  +x_T(y_A-y_D) +\alpha^2(x_Ay_D-x_Dy_A) \big).
\end{array}   \right.  \label{eq:xinB3}
\end{align}
From \eqref{eq:CirCenRad} we have that $x_T=(1-\alpha^2)x_c +\alpha^2 x_A$ and $y_T=(1-\alpha^2)y_c +\alpha^2 y_A$. Substituting these expressions into \eqref{eq:xinB3} we obtain
\begin{align}
\left.
	 \begin{array}{l l} 
 \frac{- (1\!-\!\alpha^2)[x(y_D -y_A) +y(x_A-x_D)  +y_c(x_D-x_A) + x_c(y_A-y_D) ] }{\alpha^2 (y_D-y_A) \sqrt{(x-x_D)^2+(y-y_D)^2} }  \\
 = \frac{ (1\!-\!\alpha^2)[(x_D-x_A)(y-y_c) - (y_D-y_A)(x-x_c) ] }{\alpha^2 (y_D-y_A) \sqrt{(x-x_D)^2+(y-y_D)^2} }  \\
 =0
 \end{array}   \right.  \label{eq:xinB4}
\end{align}
where equation \eqref{eq:slopes} was used in the second line of the previous equation. Thus, the choice $\omega^*=\rho$, equivalently, the aimpoint $I:(x,y)$ which is the tangent point between the orthogonal bisector and the Apollonius circle is a solution of \eqref{eq:ATapp}  and \eqref{eq:Quartic} and the Value is $V(\textbf{x}|\textbf{x}\in\mathcal{B})=0$. Furthermore, it can easily be shown that any other angle $\omega\neq \rho$ results in a terminal separation $\overline{DT}>0$. Therefore, $I:(x,y)$ is the optimal interception point and $A$ captures $T$ at the same time instant that $D$ intercepts $A$, as expected.

Now, it will be shown that the optimal solution of the ATDDG Game of Degree provides the same interception point. The optimal solution of the ATDDG Game of Degree is given by the rooting of equation (15) in Reference \cite{Garcia19TAC}. Similarly, we use the equivalent equation given by (77) in Reference \cite{Garcia19TAC}; the left hand side of such equation is given by
\begin{align}
\left.
	 \begin{array}{l l}
	 \alpha \sqrt{(x-x_T)^2+(y-y_T)^2}  \\
 -  \frac{(x-x_T)(y-y_D) - (x-x_D)(y-y_T)}{(x-x_A)(y-y_D) - (x-x_D)(y-y_A)} \sqrt{(x-x_A)^2+(y-y_A)^2} \\  
+  \frac{(x-x_T)(y-y_A) - (x-x_A)(y-y_T)}{(x-x_A)(y-y_D) - (x-x_D)(y-y_A)} \sqrt{(x-x_D)^2+(y-y_D)^2}
\end{array}   \right.   \label{eq:atddg1}
\end{align}
and we substitute $\omega^*=\rho$, equivalently, the aimpoint $I:(x,y)$ into \eqref{eq:atddg1}. Under such selection of aimpoint we have that the distance traveled by $A$ and $D$ is the same. Additionally, $\overline{TI}=\alpha\overline{AI}$. Hence, \eqref{eq:atddg1} can be written as follows
\begin{align}
\left.
	 \begin{array}{l l}
	 \frac{\sqrt{(x-x_A)^2+(y-y_A)^2}}{(x-x_A)(y-y_D) - (x-x_D)(y-y_A)}  \\
 \ \times \big( \alpha^2[(x-x_A)(y-y_D) - (x-x_D)(y-y_A)]  \\  
~~ +  (x-x_T)(y_D-y_A) - (y-y_T)(x_D-x_A)  \big)  \\
=  \frac{\sqrt{(x-x_A)^2+(y-y_A)^2}}{(x-x_A)(y-y_D) - (x-x_D)(y-y_A)}  \\
\ \times \big( (1\!-\!\alpha^2)[x(y_D -y_A) +y(x_A-x_D)] + y_T(x_D-x_A) \\
 ~~  +x_T(y_A-y_D) +\alpha^2(x_Ay_D-x_Dy_A) \big).
\end{array}   \right.   \label{eq:atddg2}
\end{align}
Since $\omega^*=\rho$ then \eqref{eq:slopes} holds. Also, the term in brackets in \eqref{eq:atddg2} is equal to the term in brackets in \eqref{eq:xinB3}, then
\begin{align}
\left.
	 \begin{array}{l l}
  \frac{\sqrt{(x-x_A)^2+(y-y_A)^2}}{(x-x_A)(y-y_D) - (x-x_D)(y-y_A)}  \\
\ \times \big( (1\!-\!\alpha^2)[x(y_D -y_A) +y(x_A-x_D)] + y_T(x_D-x_A) \\
 ~~  +x_T(y_A-y_D) +\alpha^2(x_Ay_D-x_Dy_A) \big)  \\
 = 0.
\end{array}   \right.   \nonumber
\end{align}
Similarly, the choice $\omega^*=\rho$ is a solution of equation (15) in Reference \cite{Garcia19TAC} and the Value is $V(\textbf{x}|\textbf{x}\in\mathcal{B})=0$, as expected. In conclusion, the same Value of the game is obtained by solving both Games of Degree and the optimal strategies of the Game of Degree in $\mathcal{R}_c$ and of the Game of Degree in $\mathcal{R}_e$ are the same in the Barrier surface $\mathcal{B}$.

\begin{figure}
	\begin{center}
		\includegraphics[width=8.0cm,trim=.5cm .0cm .0cm .5cm]{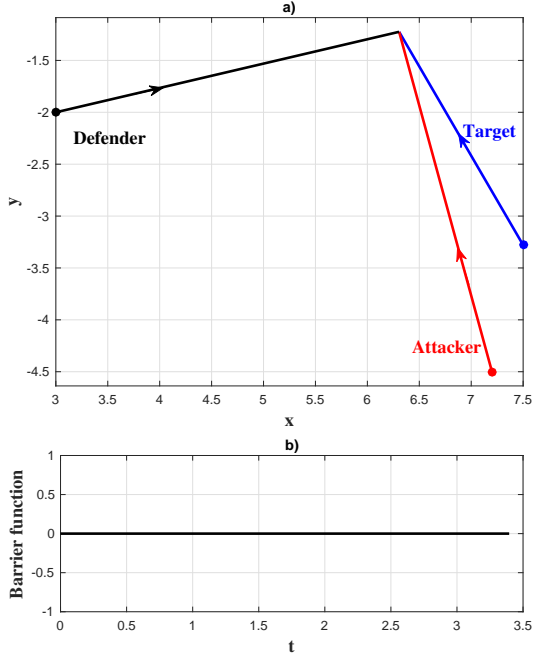}
	\caption{Example 1. a) Optimal trajectories. b) Barrier function}
	\label{fig:ExV0}
	\end{center}
\end{figure}

\section{Examples} \label{sec:examples}  
In this section we provide illustrative examples. First, an example where the state of the system is on the Barrier surface is presented. Then, the optimal strategies of the CDG are compared with respect to the strategies in Reference \cite{Liang19}.

\textit{Example 1}. Let us consider the initial positions $T=(7.5 \ -3.28)$, $A=(7.2 \ -4.5)$, $D=(3 \ -2)$, and the speed ratio $\alpha=0.7$. It holds that $\textbf{x} \in \mathcal{B}$. In this example, since $\textbf{x} \in \mathcal{B}$, the players continuously solve both Games of Degree and the optimal aimpoint is time-invariant and is the same by using either the ATDDG or the CDG: $I^*=(6.305, \ -1.224)$. The optimal trajectories are shown in Fig. \ref{fig:ExV0}.a. The Barrier function is calculated continuously for all the duration of the engagement and, as expected, $B(\textbf{x}(t);\alpha)=0$ for $0\leq t \leq t_f$ as illustrated in Fig. \ref{fig:ExV0}.b. Finally and also as one should expect, all agents meet at the interception point at the same time and the Value of the game is $V(\textbf{x})=0$.

\begin{figure}
	\begin{center}
		\includegraphics[width=7.4cm,trim=.9cm .0cm .2cm .0cm]{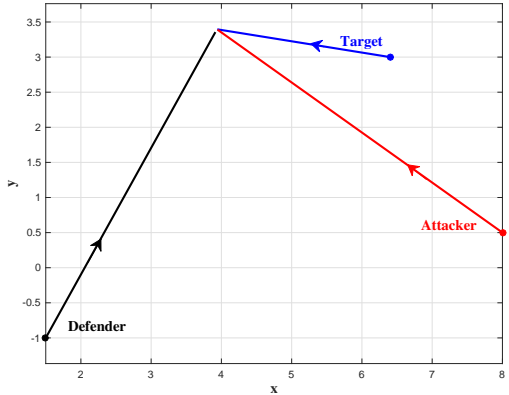}
	\caption{Example 2. Optimal play: all agents implement the optimal strategies of the CDG and the Attacker captures the Target before being intercepted by the Defender}
	\label{fig:Ex1opt}
	\end{center}
\end{figure}

\textit{Example 2}. Consider the speed ratio parameter $\alpha=0.5$ and the initial players' positions $T=(6.4 \ 3)$, $A=(8 \ 0.5)$,  and $D=(1.5 \ -1)$. According to these parameters and positions it holds that $\textbf{x} \in \mathcal{R}_c$. The agents then play the CDG and the optimal trajectories are shown in Fig. \ref{fig:Ex1opt}. Each player in this example continuously recomputes the optimal strategies \eqref{eq:OCTar}-\eqref{eq:OC} and it is obtained that the Barrier function remains positive for the duration of the engagement and that the optimal capture point is time-invariant along optimal trajectories, as expected. 

It is important to note that the Attacker, who initially was prescribed to win the game, is able to hold the state of the system in its winning region regardless of the strategies implemented by the $T/D$ team and it actually captures the Target by implementing its optimal strategy \eqref{eq:OC}. 

\begin{figure}
	\begin{center}
		\includegraphics[width=8.0cm,trim=.9cm .0cm .2cm .0cm]{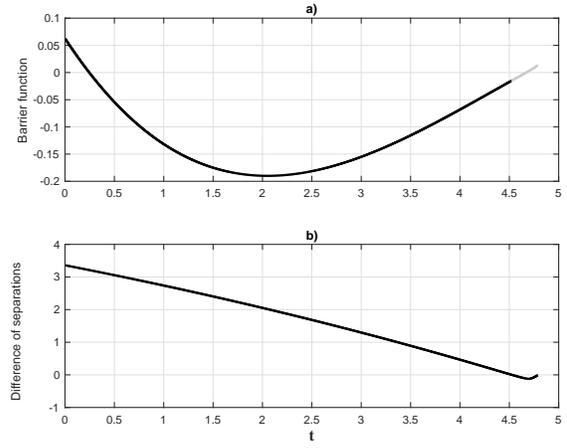}
	\caption{ Example 3. a) Barrier function. b) $\overline{AT}-\overline{DT}$}
	\label{fig:BarrSep}
	\end{center}
\end{figure}

\textit{Example 3}. We consider the same parameters and initial conditions as in Example 2. Now the Attacker implements the strategy in  \cite{Liang19},  where the seemingly optimal strategy for the Attacker is to follow Pure Pursuit (PP) on the Target. The Target and Defender will implement their cooperative optimal strategies \eqref{eq:OCTar} and \eqref{eq:OCDef} of the CDG and switch to the optimal strategies of the ATTDG when the state of the system switches regions from $\mathcal{R}_c$ to $\mathcal{R}_e$; this occurs at around $t=0.243$  as it is shown in Fig. \ref{fig:BarrSep}.a. Fig. \ref{fig:BarrSep}.b shows the difference of the distances in \eqref{eq:Red}. The system enters the subregion $\mathcal{R}_{ed}$ at around $t=4.523$. At this time $B(\textbf{x})$ becomes irrelevant and the remaining portion has been shaded in the figure ( $B(\textbf{x})$  becomes positive again at the end indicating crossing of the `$D$ branch of the hyperbola' which is the irrelevant branch as previously explained in Section \ref{subsec:gok}).  

\begin{figure}
	\begin{center}
		\includegraphics[width=8.0cm,trim=.9cm .0cm .2cm .0cm]{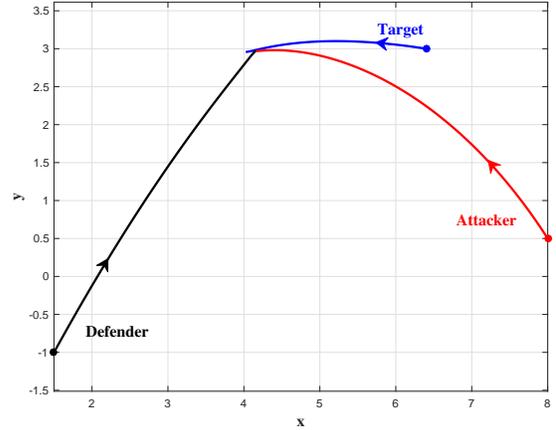}
	\caption{Example 3. Attacker implements PP and it is intercepted by the Defender before it can capture the Target}
	\label{fig:Ex1PP}
	\end{center}
\end{figure}

The trajectories when $A$ uses the PP strategy are shown in  Fig. \ref{fig:Ex1PP}. The $T$ and $D$ trajectories are not straight lines since they react to the non-optimal guidance of the Attacker. It is important to note that the Target and Defender \textit{do not know} the strategy of the Attacker; they do not know that the Attacker is implementing the PP guidance. The  $T/D$ team only uses the current state of the system, that is, the instantaneous position of the players, in order to compute the optimal strategies \eqref{eq:OCTar}-\eqref{eq:OCDef}  and defeat the Attacker.
The PP strategy may seem reasonable since it aims at minimizing capture time; however, due to the presence of the Defender we can see in Fig.  \ref{fig:Ex1PP} that $A$, who  initially was prescribed to win the game, is not able to hold the state of the system in its winning region and it ends up being intercepted by $D$ before it can capture $T$. Hence, the PP strategy obtained in Reference \cite{Liang19} does not provide a semipermeable surface, that is, the $T/D$ team is able to win the game when it was initially doomed. This is a courtesy of $A$ who erroneously implements the PP strategy. The PP strategy is \textit{not} the Attacker optimal strategy when it is pursuing the Target in the presence of the Defender.

\section{Conclusions} \label{sec:concl}  
The differential game between the Target and Defender team against the Attacker has been addressed with a focus on the Attacker's winning region. The optimal strategies of each player were synthesized and these strategies were proven to constitute the optimal solution of TAD Capture Differential Game by obtaining the Value function and showing it is continuous and continuously differentiable and it is also the solution of the Hamilton-Jacobi-Isaacs equation. The results of this paper were compared to existing approaches highlighting the superiority and optimality of the strategies obtained in this paper.


\bibliographystyle{ieeetran}
\bibliography{ReferencesTAD}

\begin{thebibliography}{10}
\providecommand{\url}[1]{#1}
\csname url@samestyle\endcsname
\providecommand{\newblock}{\relax}
\providecommand{\bibinfo}[2]{#2}
\providecommand{\BIBentrySTDinterwordspacing}{\spaceskip=0pt\relax}
\providecommand{\BIBentryALTinterwordstretchfactor}{4}
\providecommand{\BIBentryALTinterwordspacing}{\spaceskip=\fontdimen2\font plus
\BIBentryALTinterwordstretchfactor\fontdimen3\font minus
  \fontdimen4\font\relax}
\providecommand{\BIBforeignlanguage}[2]{{%
\expandafter\ifx\csname l@#1\endcsname\relax
\typeout{** WARNING: IEEEtran.bst: No hyphenation pattern has been}%
\typeout{** loaded for the language `#1'. Using the pattern for}%
\typeout{** the default language instead.}%
\else
\language=\csname l@#1\endcsname
\fi
#2}}
\providecommand{\BIBdecl}{\relax}
\BIBdecl

\bibitem{Isaacs65}
R.~Isaacs, \emph{Differential Games}.\hskip 1em plus 0.5em minus 0.4em\relax
  New York: Wiley, 1965.

\bibitem{Basar99}
T.~Basar and G.~J. Olsder, \emph{Dynamic noncooperative game theory}.\hskip 1em
  plus 0.5em minus 0.4em\relax {SIAM}, 1999, vol.~23.

\bibitem{Huang11}
H.~Huang, W.~Zhang, J.~Ding, D.~M. Stipanovic, and C.~J. Tomlin, ``Guaranteed
  decentralized pursuit-evasion in the plane with multiple pursuers,'' in
  \emph{50th IEEE Conference on Decision and Control and European Control
  Conference}, 2011, pp. 4835--4840.

\bibitem{Oyler16}
D.~W. Oyler, P.~T. Kabamba, and A.~R. Girard, ``Pursuit--evasion games in the
  presence of obstacles,'' \emph{Automatica}, vol.~65, pp. 1--11, 2016.

\bibitem{Sprinkle04}
J.~Sprinkle, J.~M. Eklund, H.~J. Kim, and S.~Sastry, ``Encoding aerial
  pursuit/evasion games with fixed wing aircraft into a nonlinear model
  predictive tracking controller,'' in \emph{43rd IEEE Conference on Decision
  and Control}, 2004, pp. 2609--2614.

\bibitem{EarlDandrea07}
M.~G. Earl and R.~D’Andrea, ``A decomposition approach to multi-vehicle
  cooperative control,'' \emph{Robotics and Autonomous Systems}, vol.~55,
  no.~4, pp. 276--291, 2007.

\bibitem{breakwell1979point}
J.~V. Breakwell and P.~Hagedorn, ``Point capture of two evaders in
  succession,'' \emph{Journal of Optimization Theory and Applications},
  vol.~27, no.~1, pp. 89--97, 1979.

\bibitem{liu2013evasion}
S.-Y. Liu, Z.~Zhou, C.~Tomlin, and K.~Hedrick, ``Evasion as a team against a
  faster pursuer,'' in \emph{American Control Conference}, 2013, pp.
  5368--5373.

\bibitem{LiCruz11}
D.~Li and J.~B. Cruz, ``Defending an asset: a linear quadratic game approach,''
  \emph{IEEE Transactions on Aerospace and Electronic Systems}, vol.~47, no.~2,
  pp. 1026--1044, 2011.

\bibitem{fisac2015pursuit}
J.~F. Fisac and S.~S. Sastry, ``The pursuit-evasion-defense differential game
  in dynamic constrained environments,'' in \emph{IEEE 54th Annual Conference
  on Decision and Control}, 2015, pp. 4549--4556.

\bibitem{Garcia18TAES}
E.~Garcia, D.~W. Casbeer, Z.~E. Fuchs, and M.~Pachter, ``Cooperative missile
  guidance for active defense of air combat vehicles,'' \emph{IEEE Transactions
  on Aerospace and Electronic Systems}, vol.~54, no.~2, pp. 706--721, 2018.

\bibitem{harini2015new}
R.~H. Venkatesan and N.~K. Sinha, ``A new guidance law for the defense missile
  of nonmaneuverable aircraft,'' \emph{IEEE Transactions on Control Systems
  Technology}, vol.~23, no.~6, pp. 2424--2431, 2015.

\bibitem{Weintraub18}
I.~Weintraub, E.~Garcia, and M.~Pachter, ``A kinematic rejoin method for active
  defense of non-maneuverable aircraft,'' in \emph{2018 American Control
  Conference}, 2018, pp. 6533--6538.

\bibitem{chen2016multiplayer}
M.~Chen, Z.~Zhou, and C.~J. Tomlin, ``Multiplayer reach-avoid games via
  pairwise outcomes,'' \emph{IEEE Transactions on Automatic Control}, vol.~62,
  no.~3, pp. 1451--1457, 2017.

\bibitem{zhou2016cooperative}
Z.~Zhou, W.~Zhang, J.~Ding, H.~Huang, D.~M. Stipanovi{\'c}, and C.~J. Tomlin,
  ``Cooperative pursuit with voronoi partitions,'' \emph{Automatica}, vol.~72,
  pp. 64--72, 2016.

\bibitem{Lorenzetti18}
J.~Lorenzetti, M.~Chen, B.~Landry, and M.~Pavone, ``Reach-avoid games via
  mixed-integer second-order cone programming,'' in \emph{57th IEEE Conference
  on Decision and Control}, 2018, pp. 4409--4416.

\bibitem{zhou2012}
Z.~Zhou, R.~Takei, H.~Huang, and C.~Tomlin, ``A general, open-loop formulation
  for reach-avoid games,'' in \emph{51st IEEE Conference on Decision and
  Control}, 2012, pp. 6501--6506.

\bibitem{margellos2011hamilton}
K.~Margellos and J.~Lygeros, ``Hamilton--jacobi formulation for reach--avoid
  differential games,'' \emph{IEEE Transactions on Automatic Control}, vol.~56,
  no.~8, pp. 1849--1861, 2011.

\bibitem{zhou2018efficient}
Z.~Zhou, J.~Ding, H.~Huang, R.~Takei, and C.~Tomlin, ``Efficient path planning
  algorithms in reach-avoid problems,'' \emph{Automatica}, vol.~89, pp. 28--36,
  2018.

\bibitem{Pachter14}
M.~Pachter, E.~Garcia, and D.~W. Casbeer, ``Active target defense differential
  game,'' in \emph{52nd Annual Allerton Conference on Communication, Control,
  and Computing}, 2014, pp. 46--53.

\bibitem{Garcia15JGCD}
E.~Garcia, D.~W. Casbeer, and M.~Pachter, ``Cooperative strategies for optimal
  aircraft defense from an attacking missile,'' \emph{Journal of Guidance,
  Control, and Dynamics}, vol.~38, no.~8, pp. 1510--1520, 2015.

\bibitem{Liang19}
L.~Liang, F.~Deng, Z.~Peng, X.~Li, and W.~Zha, ``A differential game for
  cooperative target defense,'' \emph{Automatica}, vol. 102, pp. 58--71, 2019.

\bibitem{Garcia18ACC}
{E. Garcia, D. W. Casbeer, and M. Pachter}, ``Optimal capture strategies in the
  target-attacker-defender differential game,'' in \emph{2018 American Control
  Conference}, 2018, pp. 68--73.

\bibitem{Garcia18pursuit}
E.~Garcia, D.~W. Casbeer, and M.~Pachter, ``Pursuit in the presence of a
  defender,'' \emph{Dynamic Games and Applications}, vol.~9, no.~3, pp.
  652--670, 2019.

\bibitem{Crandall84}
M.~G. Crandall, L.~C. Evans, and P.~L. Lions, ``Some properties of viscosity
  solutions of hamilton-jacobi equations,'' \emph{Transactions of the American
  Mathematical Society}, vol. 282, no.~2, pp. 487--502, 1984.

\bibitem{Lions85}
P.~L. Lions and P.~E. Souganidis, ``Differential games, optimal control and
  directional derivatives of viscosity solutions of bellman's and isaacs'
  equations,'' \emph{SIAM Journal on Control and Optimization}, vol.~23, no.~4,
  pp. 566--583, 1985.

\bibitem{getz1979qualitative}
W.~M. Getz and G.~Leitmann, ``Qualitative differential games with two
  targets,'' \emph{Journal of Mathematical Analysis and Applications}, vol.~68,
  no.~2, pp. 421--430, 1979.

\bibitem{Getz1981capturability}
W.~M. Getz and M.~Pachter, ``Capturability in a two-target game of two cars,''
  \emph{Journal of Guidance, Control, and Dynamics}, vol.~4, no.~1, pp. 15--21,
  1981.

\bibitem{ardema1985combat}
M.~D. Ardema, M.~Heymann, and N.~Rajan, ``Combat games,'' \emph{Journal of
  Optimization Theory and Applications}, vol.~46, no.~4, pp. 391--398, 1985.

\bibitem{Garcia19TAC}
E.~Garcia, D.~W. Casbeer, and M.~Pachter, ``Design and analysis of
  state-feedback optimal strategies for the differential game of active
  defense,'' \emph{IEEE Transactions on Automatic Control}, vol.~64, no.~2, pp.
  553--568, 2019.

\bibitem{Pachter16}
M.~Pachter, E.~Garcia, and D.~W. Casbeer, ``Toward a solution of the active
  target defense differential game,'' \emph{Dynamic Games and Applications},
  vol.~9, no.~1, pp. 165--216, 2019.

\bibitem{Lewin94}
J.~Lewin, \emph{Differential Games: theory and methods for solving game
  problems with singular surfaces}.\hskip 1em plus 0.5em minus 0.4em\relax
  Springer-Verlag London Limited, 1994.

\end{thebibliography}

\end{document}